\documentclass[3p,11pt,authoryear]{elsarticle}
\journal{EJOR}
\usepackage{geometry}
\newcommand{\myskip}{\smallskip}

\newcommand{\myparagraph}[1]{\paragraph{#1}}

\usepackage{appendix}

\usepackage{amsmath, amsfonts, amssymb}
\newcommand{\myper}{\! \cdot \!}
\newcommand{\myperbis}{\!\! \cdot \!\!}

\usepackage{enumerate}
\usepackage{enumitem}

\usepackage{booktabs}
\usepackage{multirow}
\usepackage{makecell}
\newcommand{\tightand}{&&{\hspace*{-.4em}}}

\usepackage[dvipsnames]{xcolor}

\usepackage{tikz}
\usetikzlibrary{
    positioning,    
    arrows.meta,    
    shapes.misc,    
    fit,            
    math,
    calc,
    backgrounds,     
    decorations.pathreplacing 
}

\usepackage[linesnumbered,ruled,vlined]{algorithm2e} 
\usepackage{setspace}

\usepackage{tcolorbox}

\usepackage[hidelinks]{hyperref}
\usepackage[nameinlink, capitalise, noabbrev]{cleveref}

\newcommand{\nbmachines}{M}
\newcommand{\nbitems}{N}
\newcommand{\nbperiods}{T}

\newcommand{\machineindex}{j}
\newcommand{\itemindex}{i}
\newcommand{\periodindex}{t}
\newcommand{\machineitemindices}{\itemindex \machineindex}
\newcommand{\machineitemindiceswithcomma}{\itemindex, \machineindex}
\newcommand{\machineitemperiodindices}{\machineitemindices \periodindex}
\newcommand{\threegenericindices}{\machineitemindices \periodindex}
\newcommand{\threegenericindiceswithcommas}{\machineitemindiceswithcomma, \periodindex}

\newcommand{\setupcost}{f}
\newcommand{\productioncost}{p}
\newcommand{\inventorycost}{h}
\newcommand{\lostsalescost}{l}
\newcommand{\demand}{d}
\newcommand{\setuptime}{s}
\newcommand{\productiontime}{b}
\newcommand{\capacity}{c}
\newcommand{\compatibility}{w}
\newcommand{\minproduction}{m}

\newcommand{\varsetup}{Y}
\newcommand{\varcarryover}{Z}
\newcommand{\varquantity}{X}
\newcommand{\varinventory}{I}
\newcommand{\varlostsales}{L}

\newcommand{\solution}{\mathcal{S}}

\newcommand{\refsuperscriptsymbol}{o}
\newcommand{\refsolution}{\solution^{\refsuperscriptsymbol}}
\newcommand{\refvarsetup}{\varsetup^{\refsuperscriptsymbol}}
\newcommand{\refvarcarryover}{\varcarryover^{\refsuperscriptsymbol}}
\newcommand{\refvarquantity}{\varquantity^{\refsuperscriptsymbol}}
\newcommand{\refvarinventory}{\varinventory^{\refsuperscriptsymbol}}

\newcommand{\maintenanceduration}{\Delta t}
\newcommand{\brokenmachineindex}{\machineindex}

\newcommand{\repairedsuperscriptsymbol}{r}
\newcommand{\repairedsolution}{\solution^{\repairedsuperscriptsymbol}}
\newcommand{\repairedvarsetup}{\varsetup^{\repairedsuperscriptsymbol}}
\newcommand{\repairedvarcarryover}{\varcarryover^{\repairedsuperscriptsymbol}}
\newcommand{\repairedvarquantity}{\varquantity^{\repairedsuperscriptsymbol}}
\newcommand{\repairedvarinventory}{\varinventory^{\repairedsuperscriptsymbol}}
\newcommand{\repairedvarlostsales}{\varlostsales^{\repairedsuperscriptsymbol}}
\newcommand{\revisedsolution}{\repairedsolution}
\newcommand{\revisedvarsetup}{\repairedvarsetup}

\usepackage{scalerel}
\newcommand{\shorttermhorizon}{\scalerel*{\tau}{t}}
\newcommand{\nbsetupchanges}{\scalerel*{\kappa}{t}}

\newcommand{\optimalsuperscriptsymbol}{\star}
\newcommand{\optimalsolution}{\solution^{\optimalsuperscriptsymbol}}
\newcommand{\optimalvarsetup}{\varsetup^{\optimalsuperscriptsymbol}}

\newcommand{\newsuperscriptsymbol}{n}
\newcommand{\newsolution}{\solution^{\newsuperscriptsymbol}}

\newcommand{\nbconvolutionblocks}{\nu}
\newcommand{\embeddingdimension}{\delta}
\newcommand{\groundtruthlabel}{y}
\newcommand{\score}{\sigma}
\newcommand{\nbselectedsetups}{\lambda}
\newcommand{\setselectedsetups}{\mathcal{\varsetup}}
\newcommand{\learningrate}{\rho}

\newcommand{\objvalue}{z}
\newcommand{\nominalobjvalue}{\objvalue^{0}}
\newcommand{\repairedobjvalue}{\objvalue^{r}}
\newcommand{\BLobjvalue}{\objvalue^{B}}
\newcommand{\GNNobjvalue}{\objvalue^{G}}
\newcommand{\bestobjvalue}{\objvalue^{\star}}
\newcommand{\Tobjvalue}{\objvalue^{T}}
\newcommand{\Pobjvalue}{\objvalue^{P}}
\newcommand{\impoverrep}{\mu}
\newcommand{\impoverrepBL}{\impoverrep(\BLobjvalue)}
\newcommand{\impoverrepGNN}{\impoverrep(\GNNobjvalue)}
\newcommand{\impoverrepLG}{\impoverrep(\bestobjvalue)}
\newcommand{\impoverrepT}{\impoverrep(\Tobjvalue)}
\newcommand{\impoverrepP}{\impoverrep(\Pobjvalue)}

\begin{document}


\begin{frontmatter}

\title{Learning to reoptimize: a GNN-aided fix-and-optimize approach\\ and an application to the Lot Sizing problem}

\author[1,3]{Mathieu Lerouge}
\author[1,2]{Andrea Lodi\corref{cor}}
\author[1]{Enrico Malaguti}
\author[1]{Michele Monaci}
\author[3]{Filippo Focacci}
\affiliation[1]{
    organization={DEI “Guglielmo Marconi” Università di Bologna}, 
    city={Bologna}, 
    country={Italy}
}
\affiliation[2]{
    organization={Jacobs Technion-Cornell Institute, Cornell Tech and Technion-IIT}, 
    city={New York}, 
    country={USA}
}
\affiliation[3]{
    organization={DecisionBrain}, 
    city={Paris},
    country={France}
}
\cortext[cor]{Corresponding author}
    
\begin{abstract}
    In many operational contexts, solutions to NP-hard combinatorial optimization problems, modeled by means of Mixed-Integer Linear Programming (MILP), may become infeasible due to unpredictable disruptions.
    Typically, reoptimizing by solving the MILP formulation on the perturbed instance is not possible as new solutions must be obtained in a very short computing time, while simple repairing heuristics may result in low-quality solutions.
    To bridge this gap, we propose a learning-to-reoptimize framework, and apply it to the Lot Sizing Problem (LSP) under machine breakdown disruptions. 
    We design a fix-and-optimize strategy aided by a Graph Neural Network (GNN) that efficiently computes a new solution within the neighborhood of a repaired solution.
    By representing the instance, the original solution and the disruption as a feature graph, we train a GNN to predict the likelihood that specific binary variables require to be modified.
    These predictions guide the selection of a small subset of variables to be reoptimized by an MILP solver, while the other variables are hard-fixed.
    Numerical experiments on a large dataset demonstrate that our approach handles effectively different problem sizes, and that it significantly outperforms a baseline alternative approach, yielding larger cost reductions within the same limited time budget.
\end{abstract}


\end{frontmatter}





\section{Introduction}

In various operational contexts, solutions to optimization problems, especially NP-hard ones, are obtained after a long computational process, lasting tens of minutes or hours. 
In many relevant settings, the temporal dimension plays a role, and decisions are made before the time at which they are actually implemented. In the time between decision and implementation, some input data may change, possibly making the decisions suboptimal or even infeasible, thus calling for recourse actions. In some cases these changes can be described by statistical distributions, allowing for the use of stochastic optimization approaches; or, at least, the decision maker can estimate intervals of variability, resulting in the possibility to adopt robust optimization schemes. 
In this paper, we consider instead situations where uncertain events are very rare but have a major impact; e.g., disruption of a connection in a transportation or telecommunication system, cancellation of a large customer demand, unavailability of a production equipment or of a major warehouse or of labor force. These situations are hardly manageable with statistical approaches and explicitly considering them at decision time could result in overconservative solutions.
\newpage

In these settings, a viable alternative is represented by reoptimization approaches, in which the decision maker computes an optimal (or near-optimal) solution to the nominal problem and possibly takes recourse actions when a disruptive and unexpected event occurs. At this time, the decision maker must compute a new high-quality solution that is feasible with respect to the disrupted instance. In many practical contexts, in order to ensure operational continuity and stability, the new solution must be close to the nominal one and be obtained in short computing time.
For this reason, a naive approach based on solving from scratch the perturbed instance is not recommendable. Consequently,  decision-makers typically resort to basic repairing approaches, such as manual adjustments or simple heuristics, to quickly restore feasibility, though these approaches generally produce low-quality solutions.

In this paper, we present a solution approach based on the idea of fixing, at the time the disruption occurs, a (possibly, large) part of the solution and exploiting the possibility to reoptimize the remaining part. 
This approach is inspired by a classical fix-and-optimize framework, an iterative approach that (i) typically requires large computing times; and (ii) may converge to a solution that is considerably different from the nominal one. To overcome these drawbacks, we use a single-iteration scheme, which turns out to be most suitable for a fast reoptimization and allows to control the deviation of the newly computed solution. The effectiveness of the resulting approach is strongly dependent on the capability to properly select the part of the solution that is fixed. We make use of Machine Learning (ML) tools for effectively performing this selection, resulting in an approach that we denote as {\em learning-to-reoptimize}.

The synergy between ML and optimization is an active and growing research field \citep{BLP21}. 
Recent works have shown that ML can effectively guide branching strategies within branch-and-bound procedures \citep{KLBS+16,LZ17} or assist in designing decompositions schemes in integer programming \citep{BCT20,XQA21}.
In the specific context of reoptimization, \cite{LMR20} and \cite{MDL24} propose using ML models to reduce the feasible solution space of the reoptimization problem, thereby accelerating solution computation in most cases.
Specifically, \cite{LMR20} predict bounds on sums of binary variables, while \cite{MDL24} predict 
whether specific binary variables should be fixed. In both cases, the considered perturbation is only on the data (capacities of facilities in \cite{LMR20} and demands of clients in \cite{MDL24}). Moreover, these approaches typically rely on ML models (such as Random Forest, Logistic Regression or Multilayer Perceptrons) that require fixed-size input vectors. 
This limits their ability to generalize across instances of varying dimensions (e.g., varying number of items).

In this paper, we make a step forward in the reoptimization direction by 1) considering more general perturbations that disrupt the structure of the problem, for example by making machines disappear, and 2) enhancing generalization through Graph Neural Networks (GNNs). 
By processing data represented as graphs, GNNs are able to capture the rich relational structure of combinatorial optimization problems and, crucially, to generalize to instances of arbitrary size.
Consequently, recent years have seen growing interest in employing GNNs within optimization contexts \citep{GCF+19,CCKL+23,ZDAZ+26}.
Notably, \cite{ZDAZ+26} use GNNs within a Deep Reinforcement Learning framework to tackle the train platforming and rescheduling problem, demonstrating the potential of graph-based learning for dynamic reoptimizing tasks.

To illustrate our methodology, we consider a specific optimization problem where time plays a crucial role, namely a Lot Sizing Problem (LSP) with setups. In detail, we consider a planning horizon divided into time periods, each characterized by the demand for a set of heterogeneous products. The problem involves a number of homogeneous machines, with setup times and penalties for unsatisfied demand.
We model this problem by means of a Mixed-Integer Linear Programming (MILP) formulation, and use a general-purpose solver for both computing an initial nominal solution and for reoptimization. 
We apply our learning-to-reoptimize framework to manage machine breakdown disruptions, where the selection of the decisions to be fixed is guided by a GNN.

The main contributions of this paper are as follows:
\begin{itemize}[nolistsep]
    \item We design a new graph representation for encoding an instance of the specific LSP application considered in the paper, along with a feasible solution and a disruption; 
    \item We propose a novel GNN-aided reoptimization strategy for fast reoptimization of a decision problem modeled as an MILP.
    By using a GNN to identify critical decision variables, we manage to reduce the practical computational effort needed to tackle the problem, thus allowing the solver to compute high-quality solutions in a short time frame;
    \item We perform a computational analysis on a large dataset, validating the generalization capabilities of our GNN and showing significant cost reduction over a baseline reoptimization strategy.
\end{itemize}

The remainder of this paper is structured as follows. 
Section \ref{sec:formulations} introduces the considered LSP and details the possible disruptions.
Section \ref{sec:reoptimization_framework} presents our proposed GNN-aided reoptimization strategy, as opposed to a baseline strategy. 
Section \ref{sec:numerical_experiments} provides a detailed computational analysis on the performance of the proposed method. 
Finally, Section \ref{sec:conclusion} concludes the paper and outlines future research directions.



\section{Mathematical formulation}
\label{sec:formulations}

In this section, we give details on the specific application under consideration and formulate it as an MILP.
Then, we introduce a specific disruption that makes a given nominal solution infeasible, and describe a straightforward deterministic procedure to derive a feasible repaired solution.
Finally, we address the reoptimization task, whose goal is to improve this repaired solution within a very short computing time.


\subsection{Nominal optimization problem}
\label{subsec:optimization_problem}

As illustrative problem, we consider a Lot Sizing application in which the decision maker is required to determine a production plan that satisfies demands over a discretized planning horizon while minimizing operational costs.
Since the founding work by \cite{WW58}, the uncapacitated single-item single-machine LSP has been extended to capacitated, multi-item, multi-machine variants, which are known to be NP-hard \citep{BY82}.
For a review on various LSP variants, the reader is referred to \cite{PW06}.
Among heuristic approaches for this class of problems, we mention the fix-and-optimize heuristic, which has been successfully applied to various LSP variants \citep{LS11,SALM13,GT18}.

In this work, we consider a multi-item, multi-machine, capacitated LSP including lost sales costs, setup costs and machine-item incompatibilities. For each setup, we require a minimum production quantity. In addition, we consider the possibility to perform a setup carry-over for (at most) one item produced on the same machine in two consecutive periods.

\myparagraph{Parameters}
An instance of our LSP consists of the following parameters:
\begin{itemize}[nolistsep, label=-]
    \item $\nbmachines, \nbitems, \nbperiods$: number of machines, items, and periods, respectively;
    \item $\setupcost_{\itemindex}$: setup cost for item $\itemindex$;
    \item $\productioncost_{\itemindex}$: unitary production cost for item $\itemindex$;
    \item $\inventorycost_{\itemindex}$: unitary inventory cost for item $\itemindex$;
    \item $\setuptime_{\itemindex}$: setup time of item $\itemindex$;
    \item $\productiontime_{\itemindex}$: unitary production time for item $\itemindex$;
    \item $\minproduction_{\itemindex}$: minimum per-setup production quantity of item $\itemindex$;
    \item $\lostsalescost_{\itemindex \periodindex}$: lost sale cost per unit of unsatisfied demand of item $\itemindex$ at period $\periodindex$;
    \item $\varinventory_{\itemindex 0}$: initial inventory of item $\itemindex$;
    \item $\demand_{\itemindex \periodindex}$: demand of item $\itemindex$ at period $\periodindex$;
    \item $\capacity_{\machineindex \periodindex}$: time availability of machine $\machineindex$ for period $\periodindex$;
    \item $\compatibility_{\itemindex \machineindex}$: compatibility between item $i$ and machine $j$ (1 if compatible, 0 otherwise).
\end{itemize}

\myparagraph{Decision variables}
The MILP formulation uses the following decision variables:
\begin{itemize}[nolistsep, label=-]
    \item $\varquantity_{\machineitemperiodindices} \geq 0$ quantity of item $\itemindex$ produced on machine $\machineindex$ at period $\periodindex$;
    \item $\varinventory_{\itemindex \periodindex} \geq 0$ inventory of item $\itemindex$ at the end of period $\periodindex$;
    \item $\varlostsales_{\itemindex \periodindex} \geq 0$ quantity of lost sales for item $\itemindex$ at the end of period $\periodindex$;
    \item $\varsetup_{\machineitemperiodindices} \in \{ 0, 1 \}$ 
    binary variable equal to 1 if machine $\machineindex$ is set up in order to produce item $\itemindex$ at period $\periodindex$, and 0 otherwise;
    \item $\varcarryover_{\machineitemperiodindices} \in \{ 0, 1 \}$ 
    binary variable equal to 1 if the production of item $\itemindex$ on machine $\machineindex$ at period $\periodindex$ is carried over 
    to period $\periodindex + 1$, and 0 otherwise.
\end{itemize}

\myparagraph{MILP formulation}
Model \eqref{eq:model_objective} -- \eqref{eq:model_non_neg_inventory_lost_sales} below defines the MILP formulation for the problem.

\vspace{-.5em}

\begin{subequations}
\renewcommand{\theequation}{\theparentequation.\arabic{equation}}
\thinmuskip=0mu
\medmuskip=0mu
\begin{align}
    \text{min} \quad &
    \sum_{\itemindex = 1}^{\nbitems} \sum_{\periodindex = 1}^{\nbperiods} \Big( \sum_{\machineindex = 1}^{\nbmachines} (\setupcost_{\itemindex} \varsetup_{\machineitemperiodindices}
    + \productioncost_{\itemindex} \varquantity_{\machineitemperiodindices})
    + \inventorycost_{\itemindex} \varinventory_{\itemindex \periodindex}
    + \lostsalescost_{\itemindex \periodindex} \varlostsales_{\itemindex \periodindex} \Big) 
    \span \span \span \span 
    \label{eq:model_objective} \\
    & \varinventory_{\itemindex (\periodindex-1)} + \sum_{\machineindex = 1}^{\nbmachines} \varquantity_{\machineitemperiodindices}
    + \varlostsales_{\itemindex \periodindex} = \demand_{\itemindex \periodindex} + \varinventory_{\itemindex \periodindex}
    \tightand \forall \; \itemindex \in 1, \; \dots, \; \nbitems, 
    \tightand \forall \; \periodindex \in 1, \; \dots, \; \nbperiods 
    \label{eq:model_balance} \\[-.6em]
    & \varlostsales_{\itemindex \periodindex} \leq \demand_{\itemindex \periodindex}
    \tightand \forall \; \itemindex \in 1, \; \dots, \; \nbitems, 
    \tightand \forall \; \periodindex \in 1, \; \dots, \; \nbperiods
    \label{eq:model_lost_sales_bound} \\
    & \sum_{\itemindex = 1}^{\nbitems} (\setuptime_{\itemindex} \varsetup_{\machineitemperiodindices}
    + \productiontime_{\itemindex} \varquantity_{\machineitemperiodindices}) \leq \capacity_{\machineindex \periodindex}
    \tightand \forall \; \machineindex \in 1, \; \dots, \; \nbmachines, 
    \tightand \forall \; \periodindex \in 1, \; \dots, \; \nbperiods 
    \label{eq:model_capacity} \\
    & \varsetup_{\machineitemperiodindices} \leq \compatibility_{\itemindex \machineindex}
    \tightand \forall \; \machineindex \in 1, \; \dots, \; \nbmachines, 
    \tightand \forall \; \itemindex \in 1, \; \dots, \; \nbitems, \quad \forall \; \periodindex \in 1, \; \dots, \; \nbperiods
    \label{eq:model_compatibility} \\
    & \varcarryover_{\machineitemperiodindices} \leq \varsetup_{\machineitemperiodindices}
    \tightand \forall \; \machineindex \in 1, \; \dots, \; \nbmachines, 
    \tightand \forall \; \itemindex \in 1, \; \dots, \; \nbitems, \quad \forall \; \periodindex \in 1, \; \dots, \; \nbperiods
    \label{eq:model_carry_over_requires_setup} \\
    & \varcarryover_{\machineitemindices (\periodindex - 1)} + \varcarryover_{\machineitemperiodindices} \leq 1
    \tightand \forall \; \machineindex \in 1, \; \dots, \; \nbmachines, 
    \tightand \forall \; \itemindex \in 1, \; \dots, \; \nbitems, \quad \forall \; \periodindex \in 1, \; \dots, \; \nbperiods
    \label{eq:model_no_consecutive_carry_over} \\
    & \varquantity_{\machineitemperiodindices} \leq M_{\threegenericindices} (\varsetup_{\machineitemperiodindices} + \varcarryover_{\machineitemindices(\periodindex - 1)})
    \tightand \forall \; \machineindex \in 1, \; \dots, \; \nbmachines, 
    \tightand \forall \; \itemindex \in 1, \; \dots, \; \nbitems, \quad \forall \; \periodindex \in 1, \; \dots, \; \nbperiods
    \label{eq:model_production_activation} \\
    & \varquantity_{\machineitemperiodindices} \geq \minproduction_{\itemindex} (\varsetup_{\machineitemperiodindices} - \varcarryover_{\machineitemperiodindices})
    \tightand \forall \; \machineindex \in 1, \; \dots, \; \nbmachines, 
    \tightand \forall \; \itemindex \in 1, \; \dots, \; \nbitems, \quad \forall \; \periodindex \in 1, \; \dots, \; \nbperiods
    \label{eq:model_min_production_when_no_carry_over} \\
    & \varquantity_{\machineitemperiodindices} + \varquantity_{\machineitemindices (\periodindex + 1)} \geq \minproduction_{\itemindex} \varcarryover_{\machineitemperiodindices}
    \tightand \forall \; \machineindex \in 1, \; \dots, \; \nbmachines, 
    \tightand \forall \; \itemindex \in 1, \; \dots, \; \nbitems, \quad \forall \; \periodindex \in 1, \; \dots, \; \nbperiods - 1
    \label{eq:model_when_carry_over} \\
    & \sum_{\itemindex = 1}^{\nbitems} \varcarryover_{\machineitemperiodindices} \leq 1
    \tightand \forall \; \machineindex \in 1, \; \dots, \; \nbmachines, 
    \tightand \forall \; \periodindex \in 1, \; \dots, \; \nbperiods
    \label{eq:model_unique_carry_over} \\
    & \varsetup_{\machineitemperiodindices}, \ \varcarryover_{\machineitemperiodindices} \in \{ 0, 1 \}
    \tightand \forall \; \machineindex \in 1, \; \dots, \; \nbmachines, 
    \tightand \forall \; \itemindex \in 1, \; \dots, \; \nbitems, \quad \forall \; \periodindex \in 1, \; \dots, \; \nbperiods
    \label{eq:model_binary} \\
    & \varcarryover_{\machineitemindices 0} = 0
    \tightand \forall \; \machineindex \in 1, \; \dots, \; \nbmachines, 
    \tightand \forall \; \itemindex \in 1, \; \dots, \; \nbitems
    \label{eq:model_zero} \\
    & \varquantity_{\machineitemperiodindices} \geq 0
    \tightand \forall \; \machineindex \in 1, \; \dots, \; \nbmachines, 
    \tightand \forall \; \itemindex \in 1, \; \dots, \; \nbitems, \quad \forall \; \periodindex \in 1, \; \dots, \; \nbperiods
    \label{eq:model_non_neg_quantity} \\
    & \varinventory_{\itemindex \periodindex}, \ \varlostsales_{\itemindex \periodindex} \geq 0
    \tightand \forall \; \itemindex \in 1, \; \dots, \; \nbitems, 
    \tightand \forall \; \periodindex \in 1, \; \dots, \; \nbperiods
    \label{eq:model_non_neg_inventory_lost_sales}
\end{align}
\label{eq:model}
\end{subequations}

The objective function \eqref{eq:model_objective} minimizes the total setup, production, inventory and lost sales costs. 
Equations \eqref{eq:model_balance} impose the inventory flow conservation, while constraints \eqref{eq:model_lost_sales_bound} define a natural upper bound on lost sales.
Constraints \eqref{eq:model_capacity} enforce machine capacity limits, whereas constraints \eqref{eq:model_compatibility} enforce machine-item compatibility.
Constraints \eqref{eq:model_carry_over_requires_setup} to \eqref{eq:model_unique_carry_over} manage the setup carry-over logic.
Specifically, constraints \eqref{eq:model_carry_over_requires_setup} ensure that a carry-over can occur at the end of a given period only after a setup is performed.
Constraints \eqref{eq:model_no_consecutive_carry_over} prevent two consecutive carry-overs for the same product on the same machine.
Constraints \eqref{eq:model_production_activation} impose that a production can be performed only when a setup or a carry-over occur. In these constraints, $M_{\threegenericindices} = \min \left( \sum_{\periodindex' = \periodindex}^{\nbperiods} \demand_{\itemindex \periodindex'}, \; (\capacity_{\machineindex \periodindex} - \setuptime_{\itemindex})/\productiontime_{\itemindex} \right)$ represents an upper bound on the production quantity.
Constraints \eqref{eq:model_min_production_when_no_carry_over} and \eqref{eq:model_when_carry_over} impose the minimum production quantity for the case without and with carry-over, respectively.
Constraints \eqref{eq:model_unique_carry_over} impose that, for each machine and period, at most one carry-over may occur.
Finally, constraints \eqref{eq:model_binary} to \eqref{eq:model_non_neg_inventory_lost_sales} define the domain of the variables, where we assume that all initial carry-over variables are fixed to 0.

\myparagraph{Nominal solution}
We assume that a high-quality feasible solution to model \eqref{eq:model_objective} -- \eqref{eq:model_non_neg_inventory_lost_sales} is available. 
Computing such a solution may require a significant effort, e.g., by running an MILP solver for a long time.
This solution, which we denote as nominal solution, corresponds to the operational plan that would be executed in the nominal case, i.e., when no disruptions occur.


\subsection{Disruption and repaired solution}
\label{subsec:disruption}

The reoptimization phase begins when an unexpected disruption invalidates the nominal solution $\refsolution$.
In this section, we first define the exact nature of the disruption that we consider in this work. We then present how to construct a repaired solution, feasible with respect to the perturbed instance, which will serve as the starting point for the reoptimization.

\myparagraph{Disruption} 
In this work, we consider two types of disruptions, both occurring at period $\periodindex = 1$ and lasting for $\maintenanceduration$ periods:
\begin{itemize}[nolistsep]
    \item \emph{Machine breakdown.}
    A single machine $\brokenmachineindex \in \{ 1, \; \dots, \; \nbmachines \}$ becomes unavailable. This disruption is modeled by setting $\capacity_{\brokenmachineindex \periodindex} = 0$ for $\periodindex = 1, \ldots, \maintenanceduration$;
    \item \emph{Plant shutdown.}
    All machines become unavailable, e.g., due to a power outage. This disruption can be modeled by setting  $\capacity_{\machineindex \periodindex}
     = 0 \quad \forall \; \machineindex \in 1, \; \dots, \; \nbmachines$ and for $\periodindex = 1, \ldots, \maintenanceduration$.
\end{itemize}
We call perturbed instance the instance affected by the disruption.
For both disruption types, the nominal solution $\refsolution$ is usually infeasible with respect to the perturbed instance. 

\myparagraph{Repaired solution}
In order to build a new solution, feasible with respect to the perturbed instance, we apply a repairing heuristic.
We call repaired solution, denoted by $\revisedsolution$, the resulting solution.
This heuristic repairing algorithm involves three main steps:
\begin{enumerate}[nolistsep, label=(\roman*)]
    \item Cancel all productions (setup, carry-over, and quantity variables) on all affected machines during  periods $\{ 1, \; \dots, \; \maintenanceduration \}$;
    \item Identify, for each machine affected by the disruption, if a carry-over has been broken and, in case, the associated product. For the resulting machine and product, redefine the production in period $\maintenanceduration + 1$ at the minimum quantity, if possible; otherwise, cancel the production;
    \item Recompute the inventory and lost sales for all affected items across the horizon.
\end{enumerate}
The repairing heuristic, detailed in \cref{alg:repairing_heuristic} in  \ref{app:repairing_heuristic}, is extremely fast in practice and returns a repaired solution $\revisedsolution$ that is feasible for the perturbed instance, though likely suboptimal.
This solution represents the starting point upon which we seek to improve, leading to the definition of a reoptimization problem, defined in the following section.

\subsection{Reoptimization problem}
\label{subsec:reoptimization_problem}

As already mentioned, the repaired solution $\revisedsolution$ provided by the heuristic algorithm described in the previous section is feasible for the perturbed instance. Better solutions can be obtained by solving formulation \eqref{eq:model_objective} -- \eqref{eq:model_non_neg_inventory_lost_sales} in the perturbed setting. However, in practical applications, organizational issues make it desirable not to deviate too much from the nominal solution planned at the beginning of the production 
period, or at least from the repaired solution (which implements minimal changes to restore feasibility). Therefore, in order to ensure operational stability, an additional constraint is imposed to force the new solution to be close to the repaired one for what concerns the setup variables in a short-term horizon right after the disruption.
More specifically, we denote by $\shorttermhorizon \leq \nbperiods$ the number of periods in the short-term horizon, and by $\nbsetupchanges \ll \nbitems \nbmachines \shorttermhorizon$ the maximum number of setup variables that can be changed.
Operational stability is enforced by imposing that the total number of setup variables within the short-term horizon that flip their value with respect to the repaired solution $\revisedsolution$ is at most $\nbsetupchanges$, and reads:
\begin{equation}
    \sum_{\itemindex = 1}^{\nbitems}
    \sum_{\machineindex = 1}^{\nbmachines}
    \Big( 
    \sum_{\substack{t = 1 \\ \text{s.t.} \\ \revisedvarsetup_{\machineitemperiodindices} = 1}}^{\shorttermhorizon} (1 - \varsetup_{\machineitemperiodindices}) +
    \sum_{\substack{t = 1 \\ \text{s.t.} \\ \revisedvarsetup_{\machineitemperiodindices} = 0}}^{\shorttermhorizon} \varsetup_{\machineitemperiodindices}
    \Big)
    \leq \nbsetupchanges.
    \label{eq:neighborhood_constraint}
\end{equation}

\myparagraph{Reoptimization problem}
In order to compute the updated production plan, we define the reoptimization problem as a new MILP obtained by adding the neighborhood constraint \eqref{eq:neighborhood_constraint} to the original formulation \eqref{eq:model_objective}--\eqref{eq:model_non_neg_inventory_lost_sales}.
The repaired solution $\revisedsolution$, though likely having a large cost, is by construction feasible for the resulting formulation, and can thus be used to warm start the solver and easy the search for an improved solution within a limited time.

The reoptimization process starting from the repaired solution forms the basis for the reoptimization strategies that will be presented in the next section and compared in the computational experiments in \cref{sec:numerical_experiments}.



\section{Reoptimization framework}
\label{sec:reoptimization_framework}

As presented in Section \ref{sec:formulations}, when a disruption occurs, it can make the high-quality nominal solution $\refsolution$ infeasible. A simple repairing heuristic may provide a new feasible solution $\repairedsolution$, but this solution is likely suboptimal. 
The goal is now to compute a new high-quality solution that improves upon $\repairedsolution$ within a short time limit.
To do this, we solve the reoptimization problem defined in \cref{subsec:reoptimization_problem}.

\myskip

In this section, we define and compare two strategies for solving this reoptimization problem, both of which use the repaired solution $\repairedsolution$ as an effective feasible warm start for the MILP solver.

\subsection{Reoptimization strategies}

\myparagraph{Baseline strategy}
The first strategy is a full-model approach.
It consists of solving the reoptimization problem (as defined in \cref{subsec:reoptimization_problem}) for a limited time.
In this approach, all the setup variables are free to be changed by the solver, bounded only by the neighborhood constraint \eqref{eq:neighborhood_constraint}.

\myparagraph{GNN-aided strategy}
The second strategy is a “fix-and-optimize” heuristic, designed to accelerate the reoptimization by drastically reducing the MILP solution space for the short-term horizon. 
The core idea is to fix the vast majority of the setup variables of the short-term horizon (i.e., periods  $\periodindex \leq \shorttermhorizon$) to their reference values in $\repairedsolution$ and to only allow a small subset of them to be free and re-optimized by the solver. 
By the temporal structure of the problem, the set-up variables in the short-term horizon are more relevant for the reoptimization process, while variables corresponding to far away periods tend to be less important and their value somehow depends on those of the former variables.
Let $\setselectedsetups$ be the subset of indices $(\threegenericindiceswithcommas)$ for setup variables $\varsetup_{\machineitemperiodindices}$ in the short-term horizon that we allow to change. We formulate a new MILP by adding the following hard-fixing constraints to the baseline model:
\begin{equation}
    \varsetup_{\machineitemperiodindices} = \repairedvarsetup_{\machineitemperiodindices} 
    \quad \forall \; \machineindex \in 1, \; \dots, \; \nbmachines, \quad \forall \; \itemindex \in 1, \; \dots, \; \nbitems, \quad \forall \; \periodindex \in 1, \; \dots, \; \shorttermhorizon, \quad (\itemindex, \machineindex, \periodindex) \notin \setselectedsetups.
    \label{eq:fixing_constraint}
\end{equation}

The feasible region of the resulting MILP corresponds to a neighborhood of the repaired solution. Though the neighborhood has an exponential size (as it happens in Very Large Neighborhood Search algorithms, see \cite{AOS00}), it can be explored successfully, provided that $\setselectedsetups$ is properly defined. 
Indeed, the challenging part of the fix-and-optimize approach consists in effectively finding a trade-off between the size of the neighborhood and the quality of the solutions it contains. Indeed, a too large $\setselectedsetups$ leaves a large number of setup variables unfixed, thus preventing an efficient exploration of the solution space, while a small subset $\setselectedsetups$ can result in solutions of poor quality. In the following section, we describe the use of ML tools to select a small subset $\setselectedsetups$ that is likely to contain the setup variables of the short-term horizon whose value should be changed to obtain a high-quality solution, while fixing all the other variables in the short-term horizon.

\myskip

\begin{figure}[htbp]
    \centering
    \colorlet{IOcolor}{Gray}
    \colorlet{commoncolor}{Gray}
    \resizebox{.7\linewidth}{!}{ 
    \begin{tikzpicture}[
        io/.style={
            rectangle, rounded corners, draw=IOcolor, fill=Gray!20,
            text centered, text width=7.5cm, minimum height=1.2cm, inner sep=3mm
        },
        inner_step/.style={
            rectangle, rounded corners,
            text centered, text width=5.5cm, minimum height=1.2cm, inner sep=3mm
        },
        common_model/.style={
            inner_step, draw=commoncolor, fill=commoncolor!20
        },
        baseline/.style={
            inner_step, draw=BurntOrange, fill=BurntOrange!20
        },
        gnn/.style={
            inner_step, draw=Green, fill=Green!20
        },
        frame/.style={
            rectangle, rounded corners, draw, dashed,
            inner sep=0.5cm,
        },
        arrow/.style={draw, -{Stealth[]}, thick}
    ]
    
        \node (reoptimization_triplet) [io] at (0, 2) {Reoptimization triplet \\ (instance, nominal solution $\refsolution$, disruption)};
        \node (reparation) [common_model] at (0, 0) {Apply repairing heuristic to $\refsolution$ \\ to obtain $\repairedsolution$}; 
        \node (common_model) [common_model] at (0, -2) {Formulate reoptimization MILP \\ \eqref{eq:model_objective} -- \eqref{eq:model_non_neg_inventory_lost_sales} and \eqref{eq:neighborhood_constraint}};
        \node (baseline_solve) [baseline] at (-4, -7) {Solve reoptimization MILP \\ using $\repairedsolution$ as warm start};
        \node (gnn) [gnn] at (4, -5) {Get GNN likelihood predictions \\ and select $\setselectedsetups$};
        \node (ml_model) [gnn] at (4, -7) {Add \eqref{eq:fixing_constraint} to the reoptimization MILP and define the reduced reopt. MILP};
        \node (ml_solve) [gnn] at (4, -9) {Solve reduced reopt. MILP \\ using $\repairedsolution$ as warm start};
        \node (end) [io] at (0, -12) {New solution $\newsolution$};
        
        \path [arrow] (reoptimization_triplet.south) -- (reparation.north);
        \path [arrow] (reparation.south) -- (common_model.north);
        \path [arrow] (common_model.south) .. controls +(1,-1) and +(0,1) .. (gnn.north);
        \path [arrow] (common_model.south) .. controls +(-1,-1) and +(0,3) .. (baseline_solve.north);
        \path [arrow] (gnn.south) -- (ml_model.north);
        \path [arrow] (ml_model.south) -- (ml_solve.north);
        \path [arrow] (baseline_solve.south) .. controls +(0,-3) and +(-1,1) .. (end.north);
        \path [arrow] (ml_solve.south) .. controls +(0,-1) and +(1,1) .. (end.north);
    
        \begin{scope}[on background layer]
            \node [frame, fit=(baseline_solve) (baseline_solve), label={[font=\bfseries, label distance=.25cm, xshift=-1.75cm]above:Baseline}] {};
            \node [frame, fit=(gnn) (ml_model) (ml_solve), label={[font=\bfseries, label distance=.25cm, xshift=1.75cm]above:GNN-aided}] {};
        \end{scope}
        
    \end{tikzpicture}
    }
    \caption{
    Comparison of baseline and GNN-aided reoptimization strategies.
    }
    \label{fig:reoptimization_strategies}
\end{figure}
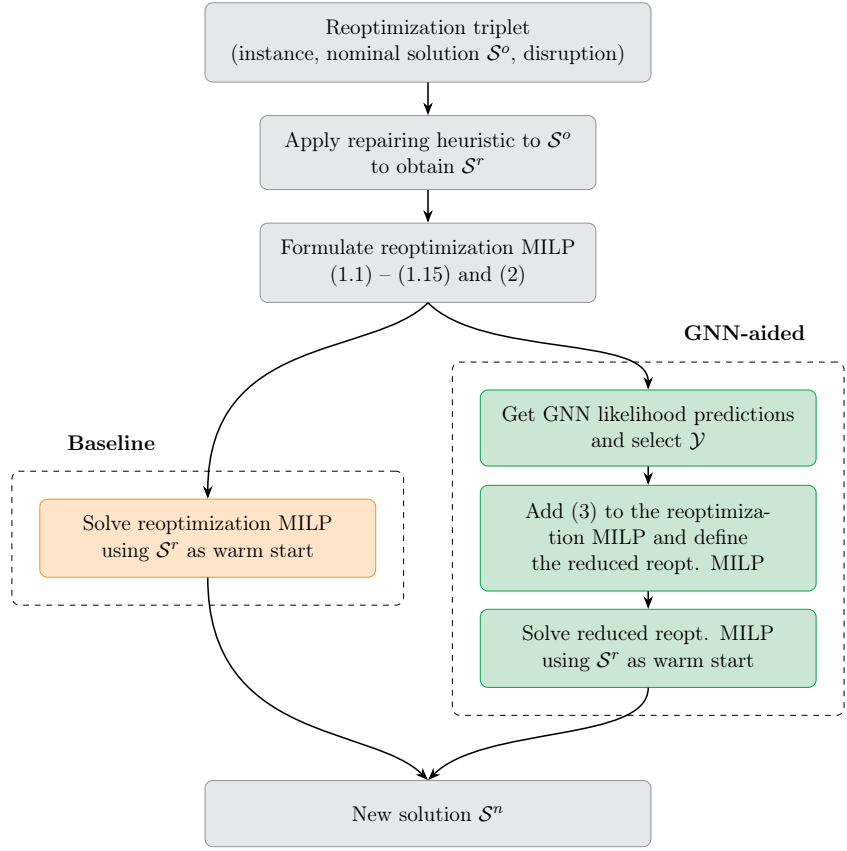

A visual comparison of the workflows of both the baseline and the GNN-aided strategies is provided in \cref{fig:reoptimization_strategies}.

\myskip

\subsection{GNN likelihood prediction}

To select the subset $\setselectedsetups$, we develop a method to predict the likelihood that each setup decision variable $\varsetup_{\machineitemperiodindices}$ of the short-term horizon ($\periodindex \leq \shorttermhorizon$) will change its value from the repaired solution $\repairedsolution$.
We design and train a GNN to perform this prediction. The GNN takes the entire reoptimization information (instance, nominal solution $\refsolution$, disruption) as a graph-structured input and outputs a change likelihood score $\score_{\machineitemperiodindices} \in [0, 1]$ for each short-term horizon variable $\varsetup_{\machineitemperiodindices}$.
The subset $\setselectedsetups$ is then formed by selecting the $\nbselectedsetups$ variables with the highest prediction scores (where $\nbselectedsetups$ is a parameter such that $\nbselectedsetups \geq \nbsetupchanges$).

\subsubsection{Graph representation}
\label{subsubsec:feature_graph}

\begin{figure}[htbp]
    \centering
    \footnotesize
    \colorlet{MPColor}{RoyalBlue}
    \colorlet{IPColor}{Green}
    \colorlet{PrColor}{BurntOrange}
    \colorlet{RsrcColor}{DarkOrchid}
    \colorlet{TimeColor}{Red}
    \colorlet{CompColor}{Salmon}
    \scalebox{.5}{
    \begin{tikzpicture}[
        base_node/.style={
            circle, draw, minimum size=1.1cm, 
            text width=1.1cm, align=center, font=\small
        },
        MP_node/.style={base_node, fill=MPColor!20, draw=MPColor, line width=1pt},
        Pr_node/.style={base_node, fill=PrColor!20, draw=PrColor, line width=1pt},
        IP_node/.style={base_node, fill=IPColor!20, draw=IPColor, line width=1pt},
        dots/.style={font={\bfseries \Large}, text=black, inner sep=0pt, minimum size=0},
        box/.style={
            rounded corners=10mm, dashed, line width=1.5pt
        },
        box_label/.style={
            font={\bfseries \Large}, text width=6cm, align=center,
        },
        arr/.style={
            draw,
            -{Stealth[]},
            thick
        },
        tw_arr/.style={
            draw,
            {Stealth[]}-{Stealth[]},
            thick
        },
        res_arr/.style={
            tw_arr,
            RsrcColor,
        },
        time_arr/.style={
            arr,
            TimeColor,
        },
        comp_mch_arr/.style={
            tw_arr,
            CompColor
        },
        comp_itm_arr/.style={
            tw_arr,
            CompColor
        }
    ]
    
        \pgfmathsetmacro{\iScale}{1.4};
        \pgfmathsetmacro{\mScale}{1.4};
        \pgfmathsetmacro{\tScale}{1.2};
        \tikzmath{
            function projectX(\i, \m, \t) {
                return (\i-1)*\iScale + (\t-1)*\tScale * 0.7; 
            };
            function projectY(\i, \m, \t) {
                return -(\m-1)*\mScale - (\t-1)*\tScale * 0.7;
            };
        };
        \pgfmathsetmacro{\tx}{\tScale * 0.7};
        \pgfmathsetmacro{\ty}{\tScale * 0.7};
    
        \begin{scope}[xshift=0, yshift=-85mm, local bounding box=MPframe]
    
            \draw [box, MPColor] ({projectX(1, 1, -.5)}, {projectY(1, 1, -.5}) rectangle ({projectX(1, 6, 7.5)}, {projectY(1, 6, 7.5});
    
            \node (mMtT) [MP_node] at ({projectX(1, 6, 6)}, {projectY(1, 6, 6)}) {$\nbmachines, \nbperiods$};
            \node (mMtt) [MP_node] at ({projectX(1, 6, 4)}, {projectY(1, 6, 4)}) {$\nbmachines, \periodindex$};
            \node (mMt2) [MP_node] at ({projectX(1, 6, 2)}, {projectY(1, 6, 2)}) {$\nbmachines, 2$};
            \node (mMt1) [MP_node] at ({projectX(1, 6, 1)}, {projectY(1, 6, 1)}) {$\nbmachines, 1$};
    
            \node (mmtT) [MP_node] at ({projectX(1, 4, 6)}, {projectY(1, 4, 6)}) {$\machineindex, \nbperiods$};
            \node (mmtt) [MP_node] at ({projectX(1, 4, 4)}, {projectY(1, 4, 4)}) {$\machineindex, \periodindex$};
            \node (mmt2) [MP_node] at ({projectX(1, 4, 2)}, {projectY(1, 4, 2)}) {$\machineindex, 2$};
            \node (mmt1) [MP_node] at ({projectX(1, 4, 1)}, {projectY(1, 4, 1)}) {$\machineindex, 1$};
            
            \node (m2tT) [MP_node] at ({projectX(1, 2, 6)}, {projectY(1, 2, 6)}) {$2, \nbperiods$};
            \node (m2tt) [MP_node] at ({projectX(1, 2, 4)}, {projectY(1, 2, 4)}) {$2, \periodindex$};
            \node (m2t2) [MP_node] at ({projectX(1, 2, 2)}, {projectY(1, 2, 2)}) {$2, 2$};
            \node (m2t1) [MP_node] at ({projectX(1, 2, 1)}, {projectY(1, 2, 1)}) {$2, 1$};
    
            \node (m1tT) [MP_node] at ({projectX(1, 1, 6)}, {projectY(1, 1, 6)}) {$1, \nbperiods$};
            \node (m1tt) [MP_node] at ({projectX(1, 1, 4)}, {projectY(1, 1, 4)}) {$1, \periodindex$};
            \node (m1t2) [MP_node] at ({projectX(1, 1, 2)}, {projectY(1, 1, 2)}) {$1, 2$};
            \node (m1t1) [MP_node] at ({projectX(1, 1, 1)}, {projectY(1, 1, 1)}) {$1, 1$};
            
            \node[dots, at={($(m2t1.south)!0.5!(mmt1.north)$)}, rotate=90] {$\dots$};
            \node[dots, at={($(mmt1.south)!0.5!(mMt1.north)$)}, rotate=90] {$\dots$}; 
            \node[dots, at={($(m2t2.south)!0.5!(mmt2.north)$)}, rotate=90] {$\dots$}; 
            \node[dots, at={($(mmt2.south)!0.5!(mMt2.north)$)}, rotate=90] {$\dots$};
            \node[dots, at={($(m2tt.south)!0.5!(mmtt.north)$)}, rotate=90] {$\dots$};
            \node[dots, at={($(mmtt.south)!0.5!(mMtt.north)$)}, rotate=90] {$\dots$}; 
            \node[dots, at={($(m2tT.south)!0.5!(mmtT.north)$)}, rotate=90] {$\dots$}; 
            \node[dots, at={($(mmtT.south)!0.5!(mMtT.north)$)}, rotate=90] {$\dots$};
            
            \node[dots, at={($(m1t2.south east)!0.5!(m1tt.north west)$)}, rotate=-45] {$\dots$};
            \node[dots, at={($(m1tt.south east)!0.5!(m1tT.north west)$)}, rotate=-45] {$\dots$}; 
            \node[dots, at={($(m2t2.south east)!0.5!(m2tt.north west)$)}, rotate=-45] {$\dots$}; 
            \node[dots, at={($(m2tt.south east)!0.5!(m2tT.north west)$)}, rotate=-45] {$\dots$};
            \node[dots, at={($(mmt2.south east)!0.5!(mmtt.north west)$)}, rotate=-45] {$\dots$};
            \node[dots, at={($(mmtt.south east)!0.5!(mmtT.north west)$)}, rotate=-45] {$\dots$}; 
            \node[dots, at={($(mMt2.south east)!0.5!(mMtt.north west)$)}, rotate=-45] {$\dots$}; 
            \node[dots, at={($(mMtt.south east)!0.5!(mMtT.north west)$)}, rotate=-45] {$\dots$};
            
        \end{scope}    
        \node [box_label, text=MPColor, yshift=7mm] at (MPframe.north) {Machine-Period nodes};
    
        \begin{scope}[xshift=85mm, yshift=0, local bounding box=IPframe]
        
            \draw [box, IPColor] ({projectX(1, 1, -.5)}, {projectY(1, 1, -.5)}) rectangle ({projectX(6, 1, 7.5)}, {projectY(6, 1, 7.5)});
            
            \node (iNtT) [IP_node] at ({projectX(6, 1, 6)}, {projectY(6, 1, 6)}) {$\nbitems, \nbperiods$};
            \node (iNtt) [IP_node] at ({projectX(6, 1, 4)}, {projectY(6, 1, 4)}) {$\nbitems, \periodindex$};
            \node (iNt2) [IP_node] at ({projectX(6, 1, 2)}, {projectY(6, 1, 2)}) {$\nbitems, 2$};
            \node (iNt1) [IP_node] at ({projectX(6, 1, 1)}, {projectY(6, 1, 1)}) {$\nbitems, 1$};
    
            \node (iitT) [IP_node] at ({projectX(4, 1, 6)}, {projectY(4, 1, 6)}) {$\itemindex, \nbperiods$};
            \node (iitt) [IP_node] at ({projectX(4, 1, 4)}, {projectY(4, 1, 4)}) {$\itemindex, \periodindex$};
            \node (iit2) [IP_node] at ({projectX(4, 1, 2)}, {projectY(4, 1, 2)}) {$\itemindex, 2$};
            \node (iit1) [IP_node] at ({projectX(4, 1, 1)}, {projectY(4, 1, 1)}) {$\itemindex, 1$};
            
            \node (i2tT) [IP_node] at ({projectX(2, 1, 6)}, {projectY(2, 1, 6)}) {$2, \nbperiods$};
            \node (i2tt) [IP_node] at ({projectX(2, 1, 4)}, {projectY(2, 1, 4)}) {$2, \periodindex$};
            \node (i2t2) [IP_node] at ({projectX(2, 1, 2)}, {projectY(2, 1, 2)}) {$2, 2$};
            \node (i2t1) [IP_node] at ({projectX(2, 1, 1)}, {projectY(2, 1, 1)}) {$2, 1$};
    
            \node (i1tT) [IP_node] at ({projectX(1, 1, 6)}, {projectY(1, 1, 6)}) {$1, \nbperiods$};
            \node (i1tt) [IP_node] at ({projectX(1, 1, 4)}, {projectY(1, 1, 4)}) {$1, \periodindex$};
            \node (i1t2) [IP_node] at ({projectX(1, 1, 2)}, {projectY(1, 1, 2)}) {$1, 2$};
            \node (i1t1) [IP_node] at ({projectX(1, 1, 1)}, {projectY(1, 1, 1)}) {$1, 1$};
            
            \node[dots, at={($(i2t1.east)!0.5!(iit1.west)$)}] {$\dots$};
            \node[dots, at={($(iit1.east)!0.5!(iNt1.west)$)}] {$\dots$}; 
            \node[dots, at={($(i2t2.east)!0.5!(iit2.west)$)}] {$\dots$}; 
            \node[dots, at={($(iit2.east)!0.5!(iNt2.west)$)}] {$\dots$};
            \node[dots, at={($(i2tt.east)!0.5!(iitt.west)$)}] {$\dots$};
            \node[dots, at={($(iitt.east)!0.5!(iNtt.west)$)}] {$\dots$}; 
            \node[dots, at={($(i2tT.east)!0.5!(iitT.west)$)}] {$\dots$}; 
            \node[dots, at={($(iitT.east)!0.5!(iNtT.west)$)}] {$\dots$};
            
            \node[dots, at={($(i1t2.south east)!0.5!(i1tt.north west)$)}, rotate=-45] {$\dots$};
            \node[dots, at={($(i1tt.south east)!0.5!(i1tT.north west)$)}, rotate=-45] {$\dots$}; 
            \node[dots, at={($(i2t2.south east)!0.5!(i2tt.north west)$)}, rotate=-45] {$\dots$}; 
            \node[dots, at={($(i2tt.south east)!0.5!(i2tT.north west)$)}, rotate=-45] {$\dots$};
            \node[dots, at={($(iit2.south east)!0.5!(iitt.north west)$)}, rotate=-45] {$\dots$};
            \node[dots, at={($(iitt.south east)!0.5!(iitT.north west)$)}, rotate=-45] {$\dots$}; 
            \node[dots, at={($(iNt2.south east)!0.5!(iNtt.north west)$)}, rotate=-45] {$\dots$}; 
            \node[dots, at={($(iNtt.south east)!0.5!(iNtT.north west)$)}, rotate=-45] {$\dots$};
            
        \end{scope}    
        \node [box_label, text=IPColor, yshift=7mm] at (IPframe.north) {Item-Period nodes};
    
        \begin{scope}[xshift=85mm, yshift=-85mm, local bounding box=PrFrame]
        
            \draw [box, PrColor] ({projectX(1, 1, -.5)}, {projectY(1, 1, -.5)}) rectangle ({projectX(6, 6, 7.5)}, {projectY(6, 6, 7.5)});
    
            \node (mMiNtT) [Pr_node] at ({projectX(6, 6, 6)}, {projectY(6, 6, 6)}) {$\nbitems, \nbmachines, \nbperiods$};
            \node (mmiNtT) [Pr_node] at ({projectX(6, 4, 6)}, {projectY(6, 4, 6)}) {$\nbitems, \machineindex, \nbperiods$};
            \node (m2iNtT) [Pr_node] at ({projectX(6, 2, 6)}, {projectY(6, 2, 6)}) {$\nbitems, 2, \nbperiods$};
            \node (mMiitT) [Pr_node] at ({projectX(4, 6, 6)}, {projectY(4, 6, 6)}) {$\itemindex, \nbmachines, \nbperiods$};
            \node (mmiitT) [Pr_node] at ({projectX(4, 4, 6)}, {projectY(4, 4, 6)}) {$\itemindex, \machineindex, \nbperiods$};
            \node (m2iitT) [Pr_node] at ({projectX(4, 2, 6)}, {projectY(4, 2, 6)}) {$\itemindex, 2, \nbperiods$};
            \node (mMi2tT) [Pr_node] at ({projectX(2, 6, 6)}, {projectY(2, 6, 6)}) {$2, \nbmachines, \nbperiods$};
            \node (mmi2tT) [Pr_node] at ({projectX(2, 4, 6)}, {projectY(2, 4, 6)}) {$2, \machineindex, \nbperiods$};
            \node (m2i2tT) [Pr_node] at ({projectX(2, 2, 6)}, {projectY(2, 2, 6)}) {$2, 2, \nbperiods$};
            
            \node (m1iitT) [Pr_node] at ({projectX(4, 1, 6)}, {projectY(4, 1, 6)}) {$\itemindex, 1, \nbperiods$};
            \node (m1iitt) [Pr_node] at ({projectX(4, 1, 4)}, {projectY(4, 1, 4)}) {$\itemindex, 1, \periodindex$};
            \node (m1iit2) [Pr_node] at ({projectX(4, 1, 2)}, {projectY(4, 1, 2)}) {$\itemindex, 1, 2$};
            \node (m1iit1) [Pr_node] at ({projectX(4, 1, 1)}, {projectY(4, 1, 1)}) {$\itemindex, 1, 1$};
            \node (m1i2tT) [Pr_node] at ({projectX(2, 1, 6)}, {projectY(2, 1, 6)}) {$2, 1, \nbperiods$};
            \node (m1i2tt) [Pr_node] at ({projectX(2, 1, 4)}, {projectY(2, 1, 4)}) {$2, 1, \periodindex$};
            \node (m1i2t2) [Pr_node] at ({projectX(2, 1, 2)}, {projectY(2, 1, 2)}) {$2, 1, 2$};
            \node (m1i2t1) [Pr_node] at ({projectX(2, 1, 1)}, {projectY(2, 1, 1)}) {$2, 1, 1$};
    
            \node (m1iNtT) [Pr_node] at ({projectX(6, 1, 6)}, {projectY(6, 1, 6)}) {$\nbitems, 1, \nbperiods$};
            \node (m1iNtt) [Pr_node] at ({projectX(6, 1, 4)}, {projectY(6, 1, 4)}) {$\nbitems, 1, \periodindex$};
            \node (m1iNt2) [Pr_node] at ({projectX(6, 1, 2)}, {projectY(6, 1, 2)}) {$\nbitems, 1, 2$};
            \node (m1iNt1) [Pr_node] at ({projectX(6, 1, 1)}, {projectY(6, 1, 1)}) {$\nbitems, 1, 1$};
            \node (m1iitT) [Pr_node] at ({projectX(4, 1, 6)}, {projectY(4, 1, 6)}) {$\itemindex, 1, \nbperiods$};
            \node (m1iitt) [Pr_node] at ({projectX(4, 1, 4)}, {projectY(4, 1, 4)}) {$\itemindex, 1, \periodindex$};
            \node (m1iit2) [Pr_node] at ({projectX(4, 1, 2)}, {projectY(4, 1, 2)}) {$\itemindex, 1, 2$};
            \node (m1iit1) [Pr_node] at ({projectX(4, 1, 1)}, {projectY(4, 1, 1)}) {$\itemindex, 1, 1$};
            \node (m1i2tT) [Pr_node] at ({projectX(2, 1, 6)}, {projectY(2, 1, 6)}) {$2, 1, \nbperiods$};
            \node (m1i2tt) [Pr_node] at ({projectX(2, 1, 4)}, {projectY(2, 1, 4)}) {$2, 1, \periodindex$};
            \node (m1i2t2) [Pr_node] at ({projectX(2, 1, 2)}, {projectY(2, 1, 2)}) {$2, 1, 2$};
            \node (m1i2t1) [Pr_node] at ({projectX(2, 1, 1)}, {projectY(2, 1, 1)}) {$2, 1, 1$};
    
            \node (mMi1tT) [Pr_node] at ({projectX(1, 6, 6)}, {projectY(1, 6, 6)}) {$1, \nbmachines, \nbperiods$};
            \node (mMi1tt) [Pr_node] at ({projectX(1, 6, 4)}, {projectY(1, 6, 4)}) {$1, \nbmachines, \periodindex$};
            \node (mMi1t2) [Pr_node] at ({projectX(1, 6, 2)}, {projectY(1, 6, 2)}) {$1, \nbmachines, 2$};
            \node (mMi1t1) [Pr_node] at ({projectX(1, 6, 1)}, {projectY(1, 6, 1)}) {$1, \nbmachines, 1$};
            \node (mmi1tT) [Pr_node] at ({projectX(1, 4, 6)}, {projectY(1, 4, 6)}) {$1, \machineindex, \nbperiods$};
            \node (mmi1tt) [Pr_node] at ({projectX(1, 4, 4)}, {projectY(1, 4, 4)}) {$\itemindex, \machineindex, \periodindex$};
            \node (mmi1t2) [Pr_node] at ({projectX(1, 4, 2)}, {projectY(1, 4, 2)}) {$1, \machineindex, 2$};
            \node (mmi1t1) [Pr_node] at ({projectX(1, 4, 1)}, {projectY(1, 4, 1)}) {$1, \machineindex, 1$};
            \node (m2i1tT) [Pr_node] at ({projectX(1, 2,  6)}, {projectY(1, 2,  6)}) {$1, 2,  \nbperiods$};
            \node (m2i1tt) [Pr_node] at ({projectX(1, 2,  4)}, {projectY(1, 2,  4)}) {$1, 2,  \periodindex$};
            \node (m2i1t2) [Pr_node] at ({projectX(1, 2,  2)}, {projectY(1, 2,  2)}) {$1, 2,  2$};
            \node (m2i1t1) [Pr_node] at ({projectX(1, 2,  1)}, {projectY(1, 2,  1)}) {$1, 2,  1$};
            
            \node (m1i1tT) [Pr_node] at ({projectX(1, 1, 6)}, {projectY(1, 1, 6)}) {$1, 1, \nbperiods$};
            \node (m1i1tt) [Pr_node] at ({projectX(1, 1, 4)}, {projectY(1, 1, 4)}) {$1, 1, \periodindex$};
            \node (m1i1t2) [Pr_node] at ({projectX(1, 1, 2)}, {projectY(1, 1, 2)}) {$1, 1, 2$};
            \node (m1i1t1) [Pr_node] at ({projectX(1, 1, 1)}, {projectY(1, 1, 1)}) {$1, 1, 1$};
    
            \node[dots, at={($(m1i2t1.east)!0.5!(m1iit1.west)$)}] {$\dots$};
            \node[dots, at={($(m1iit1.east)!0.5!(m1iNt1.west)$)}] {$\dots$};
            \node[dots, at={($(m1i2t2.east)!0.5!(m1iit2.west)$)}] {$\dots$};
            \node[dots, at={($(m1iit2.east)!0.5!(m1iNt2.west)$)}] {$\dots$};
            \node[dots, at={($(m1i2tt.east)!0.5!(m1iitt.west)$)}] {$\dots$};
            \node[dots, at={($(m1iitt.east)!0.5!(m1iNtt.west)$)}] {$\dots$};
            \node[dots, at={($(m1i2tT.east)!0.5!(m1iitT.west)$)}] {$\dots$};
            \node[dots, at={($(m1iitT.east)!0.5!(m1iNtT.west)$)}] {$\dots$};
            \node[dots, at={($(m2i2tT.east)!0.5!(m2iitT.west)$)}] {$\dots$};
            \node[dots, at={($(m2iitT.east)!0.5!(m2iNtT.west)$)}] {$\dots$};
            \node[dots, at={($(mmi2tT.east)!0.5!(mmiitT.west)$)}] {$\dots$};
            \node[dots, at={($(mmiitT.east)!0.5!(mmiNtT.west)$)}] {$\dots$};
            \node[dots, at={($(mMi2tT.east)!0.5!(mMiitT.west)$)}] {$\dots$};
            \node[dots, at={($(mMiitT.east)!0.5!(mMiNtT.west)$)}] {$\dots$};
            
            \node[dots, at={($(m2i1tT.south)!0.5!(mmi1tT.north)$)}, rotate=90] {$\dots$};
            \node[dots, at={($(mmi1tT.south)!0.5!(mMi1tT.north)$)}, rotate=90] {$\dots$};
            \node[dots, at={($(m2i1t1.south)!0.5!(mmi1t1.north)$)}, rotate=90] {$\dots$};
            \node[dots, at={($(mmi1t1.south)!0.5!(mMi1t1.north)$)}, rotate=90] {$\dots$};
            \node[dots, at={($(m2i1t2.south)!0.5!(mmi1t2.north)$)}, rotate=90] {$\dots$};
            \node[dots, at={($(mmi1t2.south)!0.5!(mMi1t2.north)$)}, rotate=90] {$\dots$};
            \node[dots, at={($(m2i1tt.south)!0.5!(mmi1tt.north)$)}, rotate=90] {$\dots$};
            \node[dots, at={($(mmi1tt.south)!0.5!(mMi1tt.north)$)}, rotate=90] {$\dots$};
            \node[dots, at={($(m2i2tT.south)!0.5!(mmi2tT.north)$)}, rotate=90] {$\dots$};
            \node[dots, at={($(mmi2tT.south)!0.5!(mMi2tT.north)$)}, rotate=90] {$\dots$};
            \node[dots, at={($(m2iitT.south)!0.5!(mmiitT.north)$)}, rotate=90] {$\dots$};
            \node[dots, at={($(mmiitT.south)!0.5!(mMiitT.north)$)}, rotate=90] {$\dots$};
            \node[dots, at={($(m2iNtT.south)!0.5!(mmiNtT.north)$)}, rotate=90] {$\dots$};
            \node[dots, at={($(mmiNtT.south)!0.5!(mMiNtT.north)$)}, rotate=90] {$\dots$};
    
            \node[dots, at={($(m1i1t2.south east)!0.5!(m1i1tt.north west)$)}, rotate=-45] {$\dots$};
            \node[dots, at={($(m1i1tt.south east)!0.5!(m1i1tT.north west)$)}, rotate=-45] {$\dots$};
            \node[dots, at={($(m1i2t2.south east)!0.5!(m1i2tt.north west)$)}, rotate=-45] {$\dots$};
            \node[dots, at={($(m1i2tt.south east)!0.5!(m1i2tT.north west)$)}, rotate=-45] {$\dots$};
            \node[dots, at={($(m1iit2.south east)!0.5!(m1iitt.north west)$)}, rotate=-45] {$\dots$};
            \node[dots, at={($(m1iitt.south east)!0.5!(m1iitT.north west)$)}, rotate=-45] {$\dots$};
            \node[dots, at={($(m1iNt2.south east)!0.5!(m1iNtt.north west)$)}, rotate=-45] {$\dots$};
            \node[dots, at={($(m1iNtt.south east)!0.5!(m1iNtT.north west)$)}, rotate=-45] {$\dots$};
            \node[dots, at={($(m2i1t2.south east)!0.5!(m2i1tt.north west)$)}, rotate=-45] {$\dots$};
            \node[dots, at={($(m2i1tt.south east)!0.5!(m2i1tT.north west)$)}, rotate=-45] {$\dots$};
            \node[dots, at={($(mmi1t2.south east)!0.5!(mmi1tt.north west)$)}, rotate=-45] {$\dots$};
            \node[dots, at={($(mmi1tt.south east)!0.5!(mmi1tT.north west)$)}, rotate=-45] {$\dots$};
            \node[dots, at={($(mMi1t2.south east)!0.5!(mMi1tt.north west)$)}, rotate=-45] {$\dots$};
            \node[dots, at={($(mMi1tt.south east)!0.5!(mMi1tT.north west)$)}, rotate=-45] {$\dots$};
    
        \end{scope}
        \node [box_label, text=PrColor, , yshift=7mm] at (PrFrame.north) {Production nodes};
    
        \begin{scope}[on background layer]
            \fill[MPColor!5] (m1t1.center) -- (m1tT.center) -- (mMtT.center) -- (mMt1.center) -- cycle;
            \fill[IPColor!5] (i1t1.center) -- (i1tT.center) -- (iNtT.center) -- (iNt1.center) -- cycle;
            \fill[MPColor!5] (m1i1t1.center) -- (m1i1tT.center) -- (mMi1tT.center) -- (mMi1t1.center) -- cycle;
            \fill[IPColor!5] (m1i1t1.center) -- (m1i1tT.center) -- (m1iNtT.center) -- (m1iNt1.center) -- cycle;
            \fill[PrColor!5] (m1i1tT) rectangle (mMiNtT);
        \end{scope}
    
        \draw[res_arr] ([xshift=1, yshift=-1]mMtT.south) to[out=-45, in=-135] ([xshift=-1, yshift=-1]mMi1tT.south);
        \draw[res_arr] ([xshift=1, yshift=-1]mMtT.south) to[out=-45, in=-135] ([xshift=-1, yshift=-1]mMi2tT.south);
        \draw[res_arr] ([xshift=1, yshift=-1]mMtT.south) to[out=-45, in=-135] ([xshift=-1, yshift=-1]$(mMi2tT.south)!0.5!(mMiitT.south)$);
        \draw[res_arr] ([xshift=1, yshift=-1]mMtT.south) to[out=-45, in=-135] ([xshift=-1, yshift=-1]mMiitT.south);
        \draw[res_arr] ([xshift=1, yshift=-1]mMtT.south) to[out=-45, in=-135] ([xshift=-1, yshift=-1]$(mMiitT.south)!0.5!(mMiNtT.south)$);
        \draw[res_arr] ([xshift=1, yshift=-1]mMtT.south) to[out=-45, in=-135] ([xshift=-1, yshift=-1]mMiNtT.south);
        \draw[res_arr] ([xshift=1, yshift=-1]iNtT.east) to[out=-45, in=45] ([xshift=1, yshift=1]m1iNtT.east);
        \draw[res_arr] ([xshift=1, yshift=-1]iNtT.east) to[out=-45, in=45] ([xshift=1, yshift=1]m2iNtT.east);
        \draw[res_arr] ([xshift=1, yshift=-1]iNtT.east) to[out=-45, in=45] ([xshift=1, yshift=1]$(m2iNtT.east)!0.5!(mmiNtT.east)$);
        \draw[res_arr] ([xshift=1, yshift=-1]iNtT.east) to[out=-45, in=45] ([xshift=1, yshift=1]mmiNtT.east);
        \draw[res_arr] ([xshift=1, yshift=-1]iNtT.east) to[out=-45, in=45] ([xshift=1, yshift=1]$(mmiNtT.east)!0.5!(mMiNtT.east)$);
        \draw[res_arr] ([xshift=1, yshift=-1]iNtT.east) to[out=-45, in=45] ([xshift=1, yshift=1]mMiNtT.east);
        \draw[time_arr] ([xshift=2]m1t1.north east) to[out=0, in=90] ([yshift=2]m1t2.north east);
        \draw[time_arr] ([xshift=2]m1t2.north east) to[out=0, in=90] ([yshift=2]$(m1t2.north east)!0.5!(m1tt.north east)$);
        \draw[time_arr] ([xshift=2]$(m1t2.north east)!0.5!(m1tt.north east)$) to[out=0, in=90] ([yshift=2]m1tt.north east);
        \draw[time_arr] ([xshift=2]m1tt.north east) to[out=0, in=90] ([yshift=2]$(m1tt.north east)!0.5!(m1tT.north east)$);
        \draw[time_arr] ([xshift=2]$(m1tt.north east)!0.5!(m1tT.north east)$) to[out=0, in=90] ([yshift=2]m1tT.north east);
        % ...
        \draw[time_arr] ([yshift=-2]i1t1.south west) to[out=-90, in=180] ([xshift=-2]i1t2.south west);
        \draw[time_arr] ([yshift=-2]i1t2.south west) to[out=-90, in=180] ([xshift=-2]$(i1t2.south west)!0.5!(i1tt.south west)$);
        \draw[time_arr] ([yshift=-2]$(i1t2.south west)!0.5!(i1tt.south west)$) to[out=-90, in=180] ([xshift=-2]i1tt.south west);
        \draw[time_arr] ([yshift=-2]i1tt.south west) to[out=-90, in=180] ([xshift=-2]$(i1tt.south west)!0.5!(i1tT.south west)$);
        \draw[time_arr] ([yshift=-2]$(i1tt.south west)!0.5!(i1tT.south west)$) to[out=-90, in=180] ([xshift=-2]i1tT.south west);
        \draw[time_arr] ([yshift=-2]mMi1t1.south west) to[out=-90, in=180] ([xshift=-2]mMi1t2.south west);
        \draw[time_arr] ([yshift=-2]mMi1t2.south west) to[out=-90, in=180] ([xshift=-2]$(mMi1t2.south west)!0.5!(mMi1tt.south west)$);
        \draw[time_arr] ([yshift=-2]$(mMi1t2.south west)!0.5!(mMi1tt.south west)$) to[out=-90, in=180] ([xshift=-2]mMi1tt.south west);
        \draw[time_arr] ([yshift=-2]mMi1tt.south west) to[out=-90, in=180] ([xshift=-2]$(mMi1tt.south west)!0.5!(mMi1tT.south west)$);
        \draw[time_arr] ([yshift=-2]$(mMi1tt.south west)!0.5!(mMi1tT.south west)$) to[out=-90, in=180] ([xshift=-2]mMi1tT.south west);
        % ...
        \draw[comp_mch_arr] ([xshift=1, yshift=1]mmiitT.east) to[out=45, in=-45] ([xshift=1, yshift=-1]m1iitT.east);
        \draw[comp_mch_arr] ([xshift=1, yshift=1]mmiitT.east) to[out=45, in=-45] ([xshift=1, yshift=-1]m2iitT.east);
        \draw[comp_mch_arr] ([xshift=1, yshift=1]mmiitT.east) to[out=45, in=-45] ([xshift=1, yshift=-1]$(m2iitT.east)!0.5!(mmiitT.east)$);
        \draw[comp_mch_arr] ([xshift=1, yshift=-1]mmiitT.east) to[out=-45, in=45] ([xshift=1, yshift=1]$(mmiitT.east)!0.5!(mMiitT.east)$);
        \draw[comp_mch_arr] ([xshift=1, yshift=-1]mmiitT.east) to[out=-45, in=45] ([xshift=1, yshift=1]mMiitT.east);
        \draw[comp_mch_arr] ([xshift=1, yshift=1]m2iitT.east) to[out=45, in=-45] ([xshift=1, yshift=-1]m1iitT.east);
        \draw[comp_mch_arr] ([xshift=1, yshift=-1]m2iitT.east) to[out=-45, in=45] ([xshift=1, yshift=1]$(m2iitT.east)!0.5!(mmiitT.east)$);
        \draw[comp_mch_arr] ([xshift=1, yshift=1]mMiitT.east) to[out=45, in=-45] ([xshift=1, yshift=-1]$(mMiitT.east)!0.5!(mmiitT.east)$);
        \draw[comp_itm_arr] ([xshift=-1, yshift=1]m1iitt.north) to[out=135, in=45] ([xshift=1, yshift=1]m1i1tt.north);
        \draw[comp_itm_arr] ([xshift=-1, yshift=1]m1iitt.north) to[out=135, in=45] ([xshift=1, yshift=1]m1i2tt.north);
        \draw[comp_itm_arr] ([xshift=-1, yshift=1]m1iitt.north) to[out=135, in=45] ([xshift=1, yshift=1]$(m1i2tt.north)!0.5!(m1iitt.north)$);
        \draw[comp_itm_arr] ([xshift=1, yshift=1]m1iitt.north) to[out=45, in=135] ([xshift=-1, yshift=1]$(m1iitt.north)!0.5!(m1iNtt.north)$);
        \draw[comp_itm_arr] ([xshift=1, yshift=1]m1iitt.north) to[out=45, in=135] ([xshift=-1, yshift=1]m1iNtt.north);
        \draw[comp_itm_arr] ([xshift=-1, yshift=1]m1i2tt.north) to[out=135, in=45] ([xshift=1, yshift=1]m1i1tt.north);
        \draw[comp_itm_arr] ([xshift=1, yshift=1]m1i2tt.north) to[out=45, in=135] ([xshift=-1, yshift=1]$(m1i2tt.north)!0.5!(m1iitt.north)$);
        \draw[comp_itm_arr] ([xshift=-1, yshift=1]m1iNtt.north) to[out=135, in=45] ([xshift=1, yshift=1]$(m1iitt.north)!0.5!(m1iNtt.north)$);
    \end{tikzpicture}
    }
    \caption{GNN feature graph structure.
    Machine-Period (resp. Item-Period and Production) nodes are represented in \textcolor{MPColor}{blue} (resp. \textcolor{IPColor}{green} and \textcolor{PrColor}{orange}).
    For readability purposes, only few arcs of each type are represented.
    Resource (resp. Time and Competition) arcs are represented in \textcolor{RsrcColor}{purple} (resp. \textcolor{TimeColor}{red} and \textcolor{CompColor}{pink}).}
    \label{fig:feature_graph}
\end{figure}
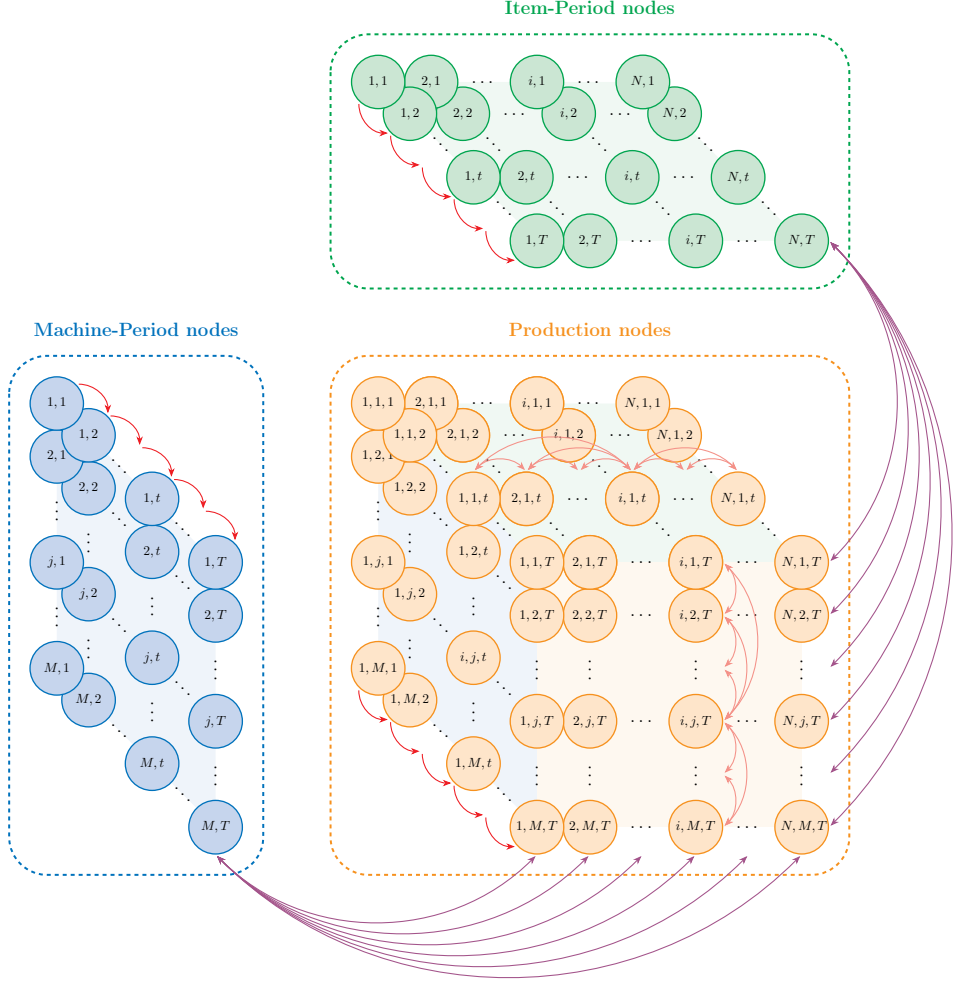

To represent the reoptimization triplet (instance, nominal solution $\refsolution$, disruption), we build a \emph{feature graph} for the GNN, as depicted in \cref{fig:feature_graph}.
The graph contains three distinct node types, namely:
\begin{itemize}[nolistsep]
    \item \emph{Machine-Period (MP) nodes}. 
    There are $\nbmachines \myper \nbperiods$ such nodes, where each node $(\machineindex, \periodindex)$ is associated with a set of features related with machine $\machineindex$ and period $\periodindex$;
    \item \emph{Item-Period (IP) nodes}. 
    There are $\nbitems \myper \nbperiods$ such nodes, where each node $(\itemindex, \periodindex)$ is associated with a set of features related with item $\itemindex$ and period $\periodindex$;
    \item \emph{Production (Pr) nodes}. 
    There are $\nbitems \myper \nbmachines \myper \nbperiods$ such nodes, where each node $(\threegenericindiceswithcommas)$ is associated with a set of features related with item $\itemindex$, machine $\machineindex$ and period $\periodindex$. 
    The GNN's prediction task is to output a score $\score_{\machineitemperiodindices}$ for each of these nodes whose period $\periodindex$ is in the short-term horizon, i.e., $\periodindex \leq \shorttermhorizon$.
\end{itemize}

The features of each node are described in \cref{tab:feature_graph_nodes}.
They encode information about 
the instance (e.g., item demand, machine capacity), 
the nominal solution (e.g., item inventory level, item production quantity), 
and the disruption (e.g., machine capacity after disruption).

\pagebreak

While most features are self-descriptive, we give specific clarification for what concerns the \emph{normalized time since end of disruption} in an MP node. 
Given a machine $\machineindex$ and a period $\periodindex$, this feature measures the time elapsed since the machine was brought back online. More into details, if a disruption lasting for $\maintenanceduration$ periods occurs on machine $\machineindex$ (i.e., starting at period 1 and ending at period $\maintenanceduration$), then this figure is computed as $(\periodindex - \maintenanceduration) / (\nbperiods - \maintenanceduration)$ for any period $\periodindex > \maintenanceduration$, and set to $0$ for any period $\periodindex \leq \maintenanceduration$.
If no disruption occurs on machine $\machineindex$, then the value is set to $-1$.

\begin{table}[htbp]
    \centering
    \renewcommand{\arraystretch}{1.2}
    \scriptsize
    \resizebox{\textwidth}{!}{
    \begin{tabular}{c c c c l c}
        \toprule
        \multicolumn{3}{c}{\textbf{Node}} &
        \multicolumn{3}{c}{\textbf{Feature}} \\
        \cmidrule(lr){1-3} \cmidrule(lr){4-6}
        \multicolumn{1}{c}{\textbf{Type}} &
        \multicolumn{1}{c}{\textbf{Indices}} &
        \multicolumn{1}{c}{\textbf{Count}} &
        \textbf{Count} &
        \textbf{Description} &
        \textbf{Value set} \\ 
        \midrule
        \multirow{8}{*}{MP} 
        & \multirow{8}{*}{$(\machineindex, \periodindex)$} 
        & \multirow{8}{*}{$\nbmachines \myper \nbperiods$} 
        & \multirow{8}{*}{8} 
        & Normalized period $\periodindex/\nbperiods$ & 
        $[0, 1]$ \\
        & & &
        & Is machine $\machineindex$ perturbed at period $\periodindex$?
        & $\{0, 1\}$ \\
        & & &
        & Normalized time since end of disruption on machine $\machineindex$ at $\periodindex$
        & $\{-1\} \cup [0, 1]$ \\
        & & &
        & Capacity of machine $\machineindex$ at period $\periodindex$ before disruption  
        & $\mathbb{R}^{+}$ \\
        & & &
        & Capacity of machine $\machineindex$ at period $\periodindex$ after disruption  
        & $\mathbb{R}^{+}$ \\
        & & &
        & Capacity left of machine $\machineindex$ at period $\periodindex$ in solution $\refsolution$  
        & $\mathbb{R}^{+}$ \\
        & & &
        & Capacity utilization of machine $\machineindex$ at period $\periodindex$ in solution $\refsolution$ 
        & $[0, 1]$ \\
        & & &
        & Total setup cost of machine $\machineindex$ at period $\periodindex$ in solution $\refsolution$  
        & $\mathbb{R}^{+}$ \\
        \midrule
        \multirow{9}{*}{IP} 
        & \multirow{9}{*}{$(\itemindex, \periodindex)$} 
        & \multirow{9}{*}{$\nbitems \myper \nbperiods$} 
        & \multirow{9}{*}{9} 
        & Normalized period $\periodindex/\nbperiods$
        & $[0, 1]$ \\
        & & & 
        & Demand $\demand_{\itemindex \periodindex}$ of item $\itemindex$ at period $\periodindex$
        & $\mathbb{R}^{+}$ \\
        & & & 
        & Minimum per-setup production quantity $\minproduction_{\itemindex}$ of item $\itemindex$;
        & $\mathbb{R}^{+}$ \\
        & & & 
        & Quantity of item $\itemindex$ produced at period $\periodindex$ in solution $\refsolution$
        & $\mathbb{R}^{+}$ \\
        & & & 
        & Inventory of item $\itemindex$ at period $\periodindex$ in solution $\refsolution$
        & $\mathbb{R}^{+}$ \\
        & & &
        & Ratio of lost sales unit cost to inventory unit cost $\lostsalescost_{\itemindex \periodindex} / \inventorycost_{\itemindex}$
        & $\mathbb{R}^{+}$ \\
        & & & 
        & Total setup cost of item $\itemindex$ at period $\periodindex$ in solution $\refsolution$
        & $\mathbb{R}^{+}$ \\
        & & & 
        & Total inventory cost of item $\itemindex$ at period $\periodindex$ in solution $\refsolution$
        & $\mathbb{R}^{+}$ \\
        & & &
        & Total lost sales cost of item $\itemindex$ at period $\periodindex$ in solution $\refsolution$
        & $\mathbb{R}^{+}$ \\
        \midrule
        \multirow{8}{*}{Pr} 
        & \multirow{8}{*}{$(\threegenericindiceswithcommas)$} 
        & \multirow{8}{*}{$\nbitems \myper \nbmachines \myper \nbperiods$} 
        & \multirow{8}{*}{8} 
        & Normalized period $\periodindex/\nbperiods$
        & $[0, 1]$ \\
        & & & 
        & Compatibility $\compatibility_{\itemindex \machineindex}$ between item $i$ and machine $j$ 
        & $\{0, 1\}$ \\
        & & & 
        & 
        Is item $\itemindex$ produced on machine $\machineindex$ at $\periodindex$ in solution $\refsolution$ ($\refvarquantity_{\machineitemperiodindices} > 0$)? \hspace{-2.5em} 
        & $\{0, 1\}$ \\
        & & & 
        & 
        Is item $\itemindex$ set on machine $\machineindex$ at period $\periodindex$ in solution $\refsolution$ ($\refvarsetup_{\machineitemperiodindices} = 1$)? 
        & $\{0, 1\}$ \\
        & & & 
        & 
        Is item $\itemindex$ carried over on $\machineindex$ from $\periodindex$ to $\periodindex \! + \! 1$ in solution $\refsolution$ ($\refvarcarryover_{\machineitemperiodindices} = 1$)? \hspace{-2.5em}
        & $\{0, 1\}$ \\
        & & & 
        & Time $\setuptime_{\itemindex} \cdot \refvarsetup_{\machineitemperiodindices}$ used to set up item $\itemindex$ on machine $\machineindex$ at $\periodindex$ in solution $\refsolution$ \hspace{-2.5em}
        & $\mathbb{R}^{+}$ \\
        & & & 
        & Quantity $\refvarquantity_{\machineitemperiodindices}$ of item $\itemindex$ produced on machine $\machineindex$ at $\periodindex$ in solution $\refsolution$ \hspace{-2.5em}
        & $\mathbb{R}^{+}$ \\
        & & & 
        & Setup cost $\setupcost_{\itemindex} \cdot \refvarsetup_{\machineitemperiodindices}$ when setting up item $\itemindex$ on $\machineindex$ at $\periodindex$ in solution $\refsolution$ \hspace{-2.5em}
        & $\mathbb{R}^{+}$ \\
        \bottomrule
    \end{tabular}
    }
    \caption{Summary of the feature graph nodes (types, feature counts, descriptions and value sets).}
    \label{tab:feature_graph_nodes}
\end{table}

The interdependencies between the information encoded in the nodes are captured by arcs, which are divided into three categories:
\begin{itemize}[nolistsep]
    \item \emph{Resources arcs}.
    These arcs connect each Pr node $(\threegenericindiceswithcommas)$ to either the MP node $(\machineindex, \periodindex)$ or to the IP node $(\itemindex, \periodindex)$.
    They allow the bidirectional flow of information between a production and its required resources (machine and item);
    \item \emph{Time arcs}.
    These arcs link nodes representing the same entity (machine or item) across consecutive periods, modeling its temporal evolution;
    \item \emph{Competition arcs}.
    These arcs connect Pr nodes that compete for the same resource (same machine and period) or contribute to the same demand (same item and period), explicitly representing production conflicts and alternatives.
\end{itemize}
The details of these connections are summarized in \cref{tab:feature_graph_edges}.

\begin{table}[htbp]
    \centering
    \renewcommand{\arraystretch}{1.2}
    \scriptsize
    \begin{tabular}{c c l c c m{6cm}}
        \toprule
        \textbf{Category} & \textbf{Source} &
        \textbf{Relation} & \textbf{Target} &
        \textbf{Count} & \textbf{Purpose} \\ 
        \midrule
        \multirow{4}{*}[-2em]{Resource} 
        & MP & \makecell[l]{involved \\ in} & Pr
        & $\nbitems \myper \nbmachines \myper \nbperiods$
        & Send information about machine $\machineindex$ at period $\periodindex$ to all production nodes that involve $\machineindex$ at $\periodindex$ \\
        & IP & \makecell[l]{involved \\ in} & Pr
        & $\nbitems \myper \nbmachines \myper \nbperiods$
        & Send information about item $\itemindex$ at period $\periodindex$ to all production nodes that involve $\itemindex$ at $\periodindex$  \\
        & Pr & involves & MP
        & $\nbitems \myper \nbmachines \myper \nbperiods$
        & Send information from all production nodes that involve machine $\machineindex$ and period $\periodindex$ to machine-period node $(\machineindex, \periodindex)$ \\
        & Pr & involves & IP
        & $\nbitems \myper \nbmachines \myper \nbperiods$
        & Send information from all production nodes that involve item $\itemindex$ and period $\periodindex$ to item-period node $(\itemindex, \periodindex)$  \\
        \midrule
        \multirow{3}{*}[-1.3em]{Time} 
        & MP & precedes & MP
        & $\nbmachines \myper (\nbperiods \! - \! 1)$
        & Propagate information about machine $\machineindex$ at period $\periodindex$ to period $\periodindex\!+\!1$ \\
        & IP & precedes & IP
        & $\nbitems \myper (\nbperiods \! - \! 1)$
        & Propagate information about item $\itemindex$ at period $\periodindex$ to period $\periodindex \! + \! 1$ \\
        & Pr & precedes & Pr
        & $\nbitems \myper \nbmachines \myper (\nbperiods \! - \! 1)$
        & Propagate the production information for item $\itemindex$ and machine $\machineindex$ at period $\periodindex$ to period $\periodindex \! + \! 1$ \\
        \midrule
        \multirow{2}{*}[-1.55em]{Competition}
        & Pr & \makecell[l]{item \\ competes \\ with} & Pr
        & $\nbitems \myper (\nbitems \! - \! 1) \myper \nbmachines \myper \nbperiods$ 
        & Send information from each production node $(\threegenericindiceswithcommas)$ to all the other production nodes involving the same machine $\machineindex$, at the same period $\periodindex$, but different items \\
        & Pr & \makecell[l]{machine \\ competes \\ with} & Pr
        & $\nbitems \myper \nbmachines \myper (\nbmachines \! - \! 1) \myper \nbperiods$ 
        & Send information from each production node $(\threegenericindiceswithcommas)$ to all the other production nodes involving the same item $\itemindex$, at the same period $\periodindex$, but different machines \\
        \bottomrule
    \end{tabular}
    \caption{Summary of the feature graph edges (types, counts and purposes).}
    \label{tab:feature_graph_edges}
\end{table}

\subsubsection{GNN architecture}
\label{sub:gnn_architecture}

\begin{figure}[htbp]
    \centering
    \scriptsize
    \resizebox{.95\linewidth}{!}{ 
    \begin{tikzpicture}[
        node distance=30mm,
        block/.style={
            rectangle,
            draw=Gray,
            fill=Gray!20,
            rounded corners=4mm,
            minimum height=22mm,
            minimum width=26mm,
            align=center,
        },
        feats/.style={
            rectangle,
            draw=BurntOrange,
            fill=BurntOrange!20,
            rounded corners=4mm,
            minimum height=22mm,
            minimum width=15mm,
            text width=22mm,
            align=right,
        },
        op/.style={ 
            rectangle,
            fill=Orchid!20,
            rounded corners=4mm,
            minimum height=38mm,
            text width=20mm,
            text depth=27mm,
            align=center
        },
        arr/.style={
            draw,
            -{Stealth[]},
            thick,
            RoyalPurple
        }
    ]
    
        \node[feats] (f1) {$\nbmachines \myperbis \nbperiods \, [8]$ \\ $\nbitems \myperbis \nbperiods \, [9]$ \\ $\nbitems \myperbis \nbmachines \myperbis \nbperiods \, [8]$};
        \node[block, anchor=east, at={([xshift=8mm]f1.west)}] (n1) {Machine-Period \\ Item-Period \\ Production};
        \node[feats, right=of f1] (f2) {$\nbmachines \myperbis \nbperiods \, [\embeddingdimension]$ \\ $\nbitems \myperbis \nbperiods \, [\embeddingdimension]$ \\ $\nbitems \myperbis \nbmachines \myperbis \nbperiods \, [\embeddingdimension]$};
        \node[block, anchor=east, at={([xshift=8mm]f2.west)}] (n2) {Machine-Period \\ Item-Period \\ Production};
        \node[feats, right=of f2] (f3) {$\nbmachines \myperbis \nbperiods \, [\embeddingdimension]$ \\ $\nbitems \myperbis \nbperiods \, [\embeddingdimension]$ \\ $\nbitems \myperbis \nbmachines \myperbis \nbperiods \, [\embeddingdimension]$};
        \node[block, anchor=east, at={([xshift=8mm]f3.west)}] (n3) {Machine-Period \\ Item-Period \\ Production};
        \node[feats, right=of f3] (f4) {$\nbitems \myperbis \nbmachines \myperbis \shorttermhorizon \, [1]$};
        \node[block, anchor=east, at={([xshift=8mm]f4.west)}] (n4) {Production};
        
        \draw[arr] (f1) -- (n2);
        \draw[arr] (f2) -- (n3);
        \draw[arr] (f3) -- (n4);
        
        \begin{scope}[on background layer]
            \node[op, anchor=south, at={($(f1.south east)!0.5!(n2.south west)$)}] (op1) {Linear \\ projection}; 
            \node[op, anchor=south, at={($(f2.south east)!0.5!(n3.south west)$)}] (op2) {Graph conv. \\ + Norm. layer \\ + Resid. co. \\ + ReLU};
            \node[op, anchor=south, at={($(f3.south east)!0.5!(n4.south west)$)}] (op3) {Filter \\ + 2-MLP \\ + Sigmoid};
        \end{scope}
        
        \draw [decorate,
               decoration={brace, amplitude=10pt, mirror, raise=8pt}, thick,
               yshift=-10pt]
            (f2.south east) -- (f3.south east)
            node [black, midway, yshift=-30pt] {Repeated $\nbconvolutionblocks$ times};
            
    \end{tikzpicture}
    }
    \caption{Schematic diagram of the GNN architecture.}
    \label{fig:GNN_architecture}
\end{figure}
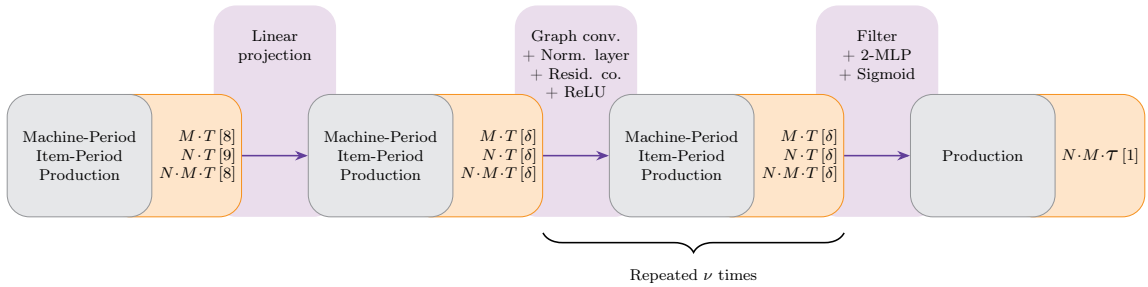

The GNN architecture is shown in \cref{fig:GNN_architecture}.
The model first projects the input features of each node, i.e., 8 features (for MP and Pr nodes) or 9 features (for IP nodes), into a common embedding space of dimension $\embeddingdimension$.
This projection is followed by a stack of $\nbconvolutionblocks$ graph convolution blocks.
Each block is heterogeneous, meaning that it learns a different message-passing function for each of the edge types (e.g., ``involved in'', ``precedes'') detailed in \cref{tab:feature_graph_edges}. 
This allows the model to learn different relational logic for resource, temporal, and competition dependencies.
Within each block, the messages passing are followed by a normalization layer and a ReLU activation.
Residual connections from the block's input to its output are also added to mitigate the vanishing gradient problem \citep{HZRS16,LMTG19}, enabling a stable training of such a deep graph network.
A final prediction head operates exclusively on the subset of Pr nodes within the short-term horizon ($\periodindex \leq \shorttermhorizon$), meaning that the model outputs prediction for exactly $\nbitems \myper \nbmachines \myper \shorttermhorizon$ nodes. 
It consists of a two-layer perceptron and a sigmoid activation function, which takes the final embedding of each node and outputs a single scalar score.
This score $\score_{\threegenericindices} \in [0, 1]$ represents the likelihood that the corresponding setup variable $\varsetup_{\threegenericindices}$ should change its value.

\myskip

By means of the scores computed by the GNN, we can determine the subset $\setselectedsetups$ used to define the reduced MILP solved in the GNN-aided strategy.
A threshold-based approach, in which we select every setup decision variable $\varsetup_{\threegenericindices}$ with a score $\score_{\threegenericindices}$ larger than the threshold, could be used. 
However, we opt for a rank-based selection to strictly control the size of subset $\setselectedsetups$.
We set a budget of $\nbselectedsetups$ variables (where $\nbselectedsetups \geq \nbsetupchanges$) and define subset $\setselectedsetups$ by selecting the $\nbselectedsetups$ variables with highest prediction scores.

\subsubsection{Ground-truth label definition}
\label{subsubsec:label_definition}

To train the GNN, a dataset of (feature graph, ground-truth label) pairs is required. The ground-truth labels must represent the ideal set of setup decision variables to change.
We assume that a dataset consisting of (instance, nominal solution, disruption) triplets is available.
For each triplet, we solve the baseline reoptimization model -- defined by equations \eqref{eq:model_objective} to \eqref{eq:model_non_neg_inventory_lost_sales} and the neighborhood constraint \eqref{eq:neighborhood_constraint} -- with a long time limit, equivalent to the time used to solve the nominal problem.
This long solving time is intended to produce a high-quality near-optimal new solution $\optimalsolution$.
We then compare this new solution $\optimalsolution$ to the repaired solution $\repairedsolution$. 
The ground-truth label $\groundtruthlabel_{\machineitemperiodindices}$ for a setup decision variable $\varsetup_{\machineitemperiodindices}$ in the short-term horizon is set to 1 if its value changed (i.e., $\optimalvarsetup_{\machineitemperiodindices} \neq \repairedvarsetup_{\machineitemperiodindices}$) and 0 otherwise. 

\myskip

This label definition process reveals a critical challenge: the resulting ground-truth labels are significantly imbalanced.
As we will see in \cref{sec:numerical_experiments}, our dataset shows that positive labels ($\groundtruthlabel_{\machineitemperiodindices} = 1$) represent only about 1\% of labels.
Such a class imbalance must be specifically addressed by our loss function.

\subsubsection{Focal loss function}
\label{subsubsec:loss_function}

Training the GNN with a standard Binary Cross-Entropy loss, $L_{BCE} = -[y \log(p) + (1-y) \log(1-p)]$, where $p$ is the model's predicted probability and $y$ is the ground-truth label, would fail on our dataset. 
The loss signal would be dominated by the large number (approximately 99\%) of “easy” negative labels, preventing the model from learning to identify the rare positive labels cases. 
To handle this significant class imbalance, we employ the Focal Loss \citep{LGG+17} defined as
\begin{equation}
    L_{FL} = -[\alpha (1-p)^\gamma y \log(p) + (1-\alpha) p^\gamma (1-y) \log(1-p)],
\end{equation}
where $\alpha \in [0, 1]$ is a class-balancing parameter and $\gamma \ge 0$ is the focusing parameter. 
The focusing parameter $\gamma$ dynamically down-weights the contribution to the loss of well-classified elements (e.g., $p \simeq 0$ for a $y=0$ label), thereby focusing training efforts on the few hard-to-classify elements (e.g., $p \simeq 0$ for a $y=1$ label), which are the ones that matter the most for our task.



\section{Numerical experiments}
\label{sec:numerical_experiments}

This section presents the numerical experiments conducted to validate the proposed GNN-aided fix-and-optimize framework.
The primary aim is to compare the solutions obtained by the GNN-aided approach against those obtained from the baseline reoptimization approach.
To this end, we first present the experimental setup, including the implementation environment and the data generation. 
We then detail the training and selection of the GNN predictor, which is a critical component of the GNN-aided strategy.
Finally, we present a comparative analysis of the two reoptimization strategies on a held-out test set.


\subsection{Experimental setup} 
\label{subsec:exprimental_setup}

In this section, we provide details about the implementation environment used for our experiments, and then describe the instance generation process and dataset preparation.

\subsubsection{Implementation environment}
\label{subsubsec:setup_environment}

All experiments were conducted on a machine equipped with an AMD Ryzen Threadripper PRO 7955WX 16-Core CPU, 64 GB of RAM, and an NVIDIA RTX A400 GPU. 
The project code was implemented in Python 3.12.3. 
The MILP models were solved using Gurobi 12.0.2 (via GurobiPy 12.0.2), while GNN models were implemented using PyTorch 2.7.0 and PyTorch Geometric 2.6.1, with CUDA 12.6 support for GPU acceleration.

\myskip

To ensure deterministic solver behavior and enhance reproducibility, we use Gurobi work units (wu)\footnotemark \ as the computational  metric for the solver runs.
\footnotetext{
Work units approximate the solver’s effort, roughly corresponding to one second of single-threaded execution according to Gurobi's documentation, though their exact value depends on hardware and model characteristics. Unlike wall-clock time, work units isolate the solver’s performance from system-level variability, such as background processes or thread contention.
}%
In our experiments, we set the short-term reoptimization budget to 10wu, while we allocate 6000wu for long-term runs to compute nominal and ground-truth solutions.

\subsubsection{Instance, solution and disruption generation}
\label{subsubsec:data_generation}

\myparagraph{Instances}
To evaluate our approach across different conditions, we generate two sets of instances, inspired by the benchmark used in \citep{Cha21}:
\begin{itemize}[nolistsep]
    \item \emph{Set 1}.
    It includes 250 instances with fixed dimensions: $\nbmachines = 3$ machines, $\nbitems = 30$ items and $\nbperiods = 30$ periods.
    For a given instance, all machines have an identical capacity $\capacity$ in $\{3000, 3500, 4000\}$. 
    No machine-item incompatibilities are considered in this set;
    \item \emph{Set 2}. 
    It includes 250 instances with varying dimensions: $\nbmachines \in \{2, 3, 4\}$ machines, $\nbitems \in \{30, 35, 40\}$ items and $\nbperiods = 30$ periods.
    Capacities are assigned as in Set 1. 
    This set introduces incompatibilities: each item is incompatible with exactly one machine.
\end{itemize}
The characteristics of these two instance sets are summarized in \cref{tab:instance_generation_main_parameters}.

\begin{table}[htbp]
    \centering
    \renewcommand{\arraystretch}{1.2}
    \scriptsize
    \begin{tabular}{l c c c c c}
        \toprule
        \multicolumn{2}{c}{\textbf{Instance}} & 
        \multicolumn{3}{c}{\textbf{Number of}} & \textbf{Machine-item} \\
        \cmidrule(lr){1-2} \cmidrule(lr){3-5}
        \textbf{Set} & \textbf{Count} & \textbf{Machines} & \textbf{Items} & \textbf{Periods} & \textbf{Incompatibilities} \\
        \midrule
        Set 1 & 250 & 3 & 30 & 30 & None \\
        Set 2 & 250 & $\{2, 3, 4\}$ & $\{30, 35, 40\}$ & 30 & 1 per item \\
        \bottomrule
    \end{tabular}
    \caption{Characteristics of the generated instance sets.}
    \label{tab:instance_generation_main_parameters}
\end{table}

For both sets, items are categorized into three groups based on their demands and costs:
\begin{itemize}[nolistsep]
    \item \emph{High priority} (6 to 8 items). 
    The total demand per period of this group accounts for 40--50\% of the total machine capacity per period. 
    The items have low inventory unit costs in $[0.05, 0.15]$, but very high initial lost sales unit costs in $[5, 9]$, decreasing over time to $[0.1, 0.2]$;
    \item \emph{Medium priority} (9 to 11 items).
    The total demand per period of this group accounts for 40--50\% of the total machine capacity per period. 
    The items have higher inventory costs in $[0.2, 0.3]$, and moderate initial lost sales costs in $[0.9, 1.1]$, decreasing to $[0.1, 0.2]$;
    \item \emph{Low priority} (remaining items). 
    The total demand per period of this group of items accounts for 20--30\% of the total machine capacity per period.
    These items have low inventory costs in $[0.05, 0.35]$, and moderate initial lost sales costs in $[0.9, 1.1]$, decreasing to $[0.1, 0.2]$.
\end{itemize}
Unit lost sales costs are decreasing over time to reflect the different impact of a lost sale in a forthcoming period (which is likely to result in a real loss) with respect to a lost sale in a far-away period (for which the sale may refer to a forecast demand).
By construction, the cumulative demand of these three priority groups is designed to slightly exceed the total available machine capacity, leading to difficult trade-offs between production and lost sales.
The characteristics of these item categories are summarized in \cref{tab:instance_generation_item_categories}.

\begin{table}[htbp]
    \centering
    \renewcommand{\arraystretch}{1.2}
    \scriptsize
    \begin{tabular}{l c c c c}
        \toprule
        \textbf{Priority} & \textbf{Item count} & \textbf{Capacity share} & \textbf{Unit inv. cost} & \textbf{Unit lost sales cost} \\
        \midrule
        High & $6-8$ & $40\% - 50 \%$ & $[0.05, 0.15]$ & $[5, 9] \to [0.1, 0.2]$ \\
        Medium & $9-11$ & $40\% - 50\%$ & $[0.20, 0.30]$ & $[0.9, 1.1] \to [0.1, 0.2]$ \\
        Low & remainder & $20\% - 30\%$ & $[0.05, 0.35]$ & $[0.9, 1.1] \to [0.1, 0.2]$ \\
        \bottomrule
    \end{tabular}
    \caption{Characteristics of the item categories.}
    \label{tab:instance_generation_item_categories}
\end{table}

\myskip

Finally, to ensure diverse and challenging scenarios, other parameters are generated as follows:
\begin{itemize}[nolistsep]
    \item \emph{Demands}. 
    Item demands are time-varying.
    They are generated independently for each period to fluctuate around the target capacity share defined by the priority category;
    \item \emph{Production}. 
    Without loss of generality, unit production times are set to 1 and unit production costs to 0;
    \item \emph{Setups}. 
    Setup times are generated to consume between 10\% and 20\% of a machine's capacity per period. Setup costs are scaled to be approximately 10\% of the corresponding setup time;
    \item \emph{Minimum Production}. 
    The minimum production quantity is randomized within a range of $[60\%, 140\%]$ of the item's average demand, ensuring that production batches are significant relative to demand.
\end{itemize}

\myparagraph{Nominal solutions}
For each of the 500 generated instances, we solve the nominal problem \eqref{eq:model_objective} -- \eqref{eq:model_non_neg_inventory_lost_sales} using a long computing budget of 6000wu.
This provides a high-quality nominal solution $\refsolution$ for each instance, which would serve as the operational plan to execute, assuming no disruption occurs.

\myparagraph{Disruptions}
We generate multiple disruption scenarios for each (instance, solution) pair. 
We consider two types of disruptions:
\begin{itemize}[nolistsep]
    \item \emph{Machine Breakdown (MB)}.
    A single machine becomes unavailable for a duration $\maintenanceduration$ of 4 or 5 periods;
    \item \emph{Plant Shutdown (PS).} 
    All machines become unavailable for a duration $\maintenanceduration$ of 1 or 2 periods.
\end{itemize}
These disruptions make the nominal solution infeasible, calling for a reoptimization process.

\myskip
\myskip

To quantify the impact of these disruptions, \cref{tab:disruption_impact} summarizes the deterioration in solution quality.
Let $\nominalobjvalue$ and $\repairedobjvalue$ denote the objective values of the nominal solution $\refsolution$ and the repaired solution $\revisedsolution$ computed by Algorithm \ref{alg:repairing_heuristic} (see Section \ref{subsec:reoptimization_problem} and Appendix \ref{app:repairing_heuristic}), respectively.
We report the cost increase of the repaired solution with respect to the nominal one, computed as $(\repairedobjvalue - \nominalobjvalue)/\nominalobjvalue$. 
To measure the operational instability, the table also gives the average number of setup variables that differ between $\refsolution$ and $\revisedsolution$. This number is given both as an absolute count and as a percentage relative to the total number of setups in $\refsolution$.

\begin{table}[htbp]
    \centering
    \renewcommand{\arraystretch}{1.2}
    \scriptsize
    \begin{tabular}{l r r r}
        \toprule
        \textbf{Instance} & 
        \multicolumn{1}{l}{\textbf{Cost}} & 
        \multicolumn{2}{c}{\textbf{Setup modifications}} \\
        \cmidrule(lr){3-4}
        \multicolumn{1}{l}{\textbf{Set}} & 
        \multicolumn{1}{l}{\textbf{Increase}} & 
        \textbf{Count} & \textbf{Relative} \\
        \midrule
        Set 1 & 
        104.18\% & 7.94 & 4.17\% \\
        Set 2 & 
        98.53\% & 8.34 & 4.38\% \\
        \midrule
        All & 
        101.45\% & 8.13 & 4.27\% \\
        \bottomrule
    \end{tabular}
    \caption{Disruption impact on cost increase and setups relative to nominal solution.}
    \label{tab:disruption_impact}
\end{table}

The results in \cref{tab:disruption_impact} reveal a critical insight: while the repairing heuristic modifies a limited number of setups ($4.27\%$ of the setups in $\refsolution$ on average), the solution worsening is significant, with costs doubling on average ($+101.45\%$).
This relevant degradation is typically due to high penalties associated with lost sales that the repairing heuristic cannot efficiently mitigate.
Consequently, the repaired solution $\revisedsolution$ is unsatisfactory, necessitating a more sophisticated reoptimization process to recover solution quality.

\subsubsection{Preparing datasets}
\label{subsubsec:datasets}

The reoptimization problem definition includes the neighborhood constraint \eqref{eq:neighborhood_constraint}, which depends on two parameters $\shorttermhorizon$ and $\nbsetupchanges$.
To simplify the analysis and interpretation of the results, in all our instances, we used the same values for these parameters and set the short-term horizon to $\shorttermhorizon = 10$ periods and the maximum number of setup changes to $\nbsetupchanges = 10$. This implies that a very small fraction of the setup variables in the first $\shorttermhorizon$ periods are allowed to change; e.g., for instances in Set 1, the first $\shorttermhorizon$ periods involve 900 setup variables, among which only 10 (i.e., around 1\%) are allowed to change their value. Though the resulting neighborhood constraint is binding, the reoptimization problem includes feasible solutions whose value considerably improves over the repaired solution $\repairedsolution$, as shown in \cref{subsec:comparison}.

\myskip

For each of the (instance, nominal solution, disruption) triplets that we have generated, we compute the corresponding feature graph as detailed in \cref{sec:reoptimization_framework}.
These graphs serve as inputs for our GNN predictor.

To ensure a rigorous evaluation, we partition all the data at the instance level.
Each set of instances (Set 1 and Set 2) is split into two subsets A and B, as follows.
\begin{itemize}[nolistsep]
    \item \emph{Subset A.} 
    This subset includes 150 instances from Set 1 and 150 instances from Set 2, and is used for the selection of the GNN configuration.
    All the (instance, nominal solution, disruption) triplets generated from these instances are used exclusively for training (70\%), validating (15\%) and testing (15\%) the GNN;
    \item \emph{Subset B.} 
    This held-out subset includes 100 instances from Set 1 and 100 instances from Set 2, and is reserved exclusively for the final comparative analysis of reoptimization strategies, presented in \cref{subsec:comparison}.
    The GNN never encounters the related data during the model selection phase.
\end{itemize}

To ensure that both learning and evaluation of the final GNN  are performed across uniform and variable instance structure, the selection is done using the union of subsets A from both Set 1 and Set 2. Conversely, the final comparative analysis is performed on the union of subsets B only.


\subsection{Selection of GNN configuration}
\label{subsec:GNN_selection}

This section describes the process of training, testing and selecting the final GNN configuration used in the GNN-aided strategy.
This is performed entirely on reoptimization triplets related to instances from subset A.

\subsubsection{Ground-truth label computation}
\label{subsubsec:label_computation}

As detailed in \cref{subsubsec:label_definition}, the ground-truth labels are computed by solving the reoptimization problem with a long computing budget (6000wu) to find a near-optimal solution $\optimalsolution$.
For each setup decision variable $\varsetup_{\machineitemperiodindices}$ of the short-term horizon ($\periodindex \leq \shorttermhorizon$), the ground-truth label $\groundtruthlabel_{\machineitemperiodindices}$ is set to 1 if $\varsetup_{\machineitemperiodindices}$ changes between the repaired solution $\repairedsolution$ and the near-optimal one $\optimalsolution$ (i.e., $\optimalvarsetup_{\machineitemperiodindices} \neq \repairedvarsetup_{\machineitemperiodindices}$), and is set to 0 otherwise.
This computation is performed for all triples from instances in subset A.

\myskip

As anticipated, the resulting labels are extremely imbalanced: across our training set, only 1.05\% of the ground-truth labels within the short-term horizon are positive (i.e., equal to 1).
This confirms the necessity of the focal loss function discussed in \cref{subsubsec:loss_function}.

\subsubsection{Model training, validation and selection}
\label{subsubsec:GNN_selection}

We considered several GNN and training configurations by varying hyperparameters, as summarized in \cref{tab:hyperparameters}.
The parameters we tuned include the embedding dimension $\embeddingdimension$, the number of graph convolution blocks $\nbconvolutionblocks$, the initial learning rate $\learningrate$, and the focal loss parameters $\alpha$ and $\gamma$.

\myskip

\begin{table}[htbp]
    \centering
    \renewcommand{\arraystretch}{1.2}
    \scriptsize
    \begin{tabular}{c l l}
        \toprule
        \multicolumn{2}{l}{\textbf{Hyperparameter}} & \textbf{Tested values} \\
        \midrule
        $\embeddingdimension$ 
        & GNN embedding dimension 
        & $\{ 32, 64 \}$ \\
        $\nbconvolutionblocks$ 
        & GNN number of convolution blocks 
        & $\{ 3, 4 \}$ \\
        $\learningrate$ 
        & Adam optimizer initial learning rate 
        & $\{ 0.0005, 0.0002 \}$ \\
        $\alpha$ 
        & Focal loss class-balancing parameter 
        & $\{ 0.25, 0.5 \}$ \\
        $\gamma$ 
        & Focal loss focusing parameter 
        & $\{ 2, 3 \}$ \\
        \bottomrule
    \end{tabular}
    \caption{Hyperparameters and the values tested in the grid search.}
    \label{tab:hyperparameters}
\end{table}

\myskip

We evaluated the performance of the GNN for different configurations of the parameters, by using standard classification metrics based on the confusion matrix components:
\begin{itemize}[nolistsep, label=-]
    \item \emph{true positives} ($TP$), i.e., number of setups that are correctly predicted to change;
    \item \emph{true negatives} ($TN$), i.e., number of setups that are correctly predicted to remain fixed;
    \item \emph{false positives} ($FP$), i.e., number of setups that are predicted to change but should have remained fixed;
    \item \emph{false negatives} ($FN$), i.e., number of setups that are predicted to remain fixed but should have changed.
\end{itemize}

\medskip

A commonly used metric for evaluating the performance of a GNN is the {\em accuracy}, defined as
\begin{equation*}
    \text{accuracy} = \frac{TP + TN}{TP + TN + FP + FN}.
\end{equation*}
Accuracy represents a measure of how often the GNN correctly predicts the labels. However, given the significant class imbalance previously described (where positive labels are about 1\%), the accuracy tends to be a meaningless metric. 
Indeed, a trivial model predicting no change for setup variables would achieve near-perfect accuracy (about 99\%) but fail to identify any relevant setup variable fixes.

For this reason, we focus on the following metrics
\begin{equation*}
    \text{precision} = \frac{TP}{TP + FP}, 
    \quad
    \text{recall} = \frac{TP}{TP + FN}.
\end{equation*}

The values of these metrics describe the behavior of our GNN network and determine the performance of the ``fix-and-optimize'' heuristic.
\begin{itemize}[nolistsep]
    \item {\em Precision} indicates how reliable the model’s positive predictions are.
    A high precision corresponds to a low value for $FP$. This indicates that subset $\setselectedsetups$ does not contain too many irrelevant setup variables, which would otherwise increase the solution space of the reoptimization problem and prevent the solver from finding high-quality solutions in a limited time.    
    \item {\em Recall} reports how well the model is able to identify all the relevant (positive) instances.
    A high recall corresponds to a low value for $FN$. This indicates that the setup variables that need to change are included in the subset of selected free setup variables $\setselectedsetups$, thus allowing the MILP solver to identify setup changes to apply and find high-quality solutions.
\end{itemize}

In addition, we make use of the F1-score, defined as
\begin{equation*}
    \text{F1-score} = 2 \ \frac{\text{precision} \cdot \text{recall}}{\text{precision} + \text{recall}}.
\end{equation*}
This value is the harmonic mean of precision and recall, providing a balanced aggregate measure.
While we report the F1-score, our selection policy is based on the \emph{highest recall} achieved on the validation set, subject to a \emph{precision of at least 33\%} (ensuring that at least one out of three predicted changes is relevant).

\myskip

We performed a grid search over the possible configurations in \cref{tab:hyperparameters} by training each model on the training set of subset A and evaluating its performance on the validation set of subset A.
The configuration that achieved the best score on our primary selection metric was selected as the final configuration, and corresponds to $\embeddingdimension = 64$, $\nbconvolutionblocks = 4$, $\learningrate = 0.0005$, $\alpha = 0.1$, and $\gamma = 2$.
On the test set of subset A, this configuration achieved a final F1-score of $45.88\%$, with a recall of $75.71\%$ and a precision of $32.91\%$.
This configuration has been used for all the remaining experiments. It is worth noting that these results corroborate the choice of using a GNN architecture for our prediction problem. Indeed, we have worked on an heterogeneous dataset with instances of different sizes (number of items, machines, disruptions) in a rather transparent way, without noticing limitations in our evaluation metrics.


\subsection{Comparative analysis of the reoptimization strategies}
\label{subsec:comparison}

We now present the comparison of the two reoptimization strategies presented in \cref{sec:reoptimization_framework}.
This experiment is performed on subset B, i.e., the held-out set of instances whose related reoptimization triplets have never been used during the GNN selection phase.
Both the baseline and the GNN-aided strategies are run on each triplet related to subset B with an identical short time frame of 10wu.
For the GNN-aided strategy, we use the best GNN configuration selected in \cref{subsubsec:GNN_selection}, and a setup selection budget that we fix to $\nbselectedsetups = 3 \nbsetupchanges = 30$ to account for the observed precision of the GNN.

\myskip

Our computational analysis is based on the comparison with two reference solutions:
(i) the best-known solution, obtained from the long-term run, having value $\bestobjvalue$; and (ii) the repaired solution having value $\repairedobjvalue$.
For a given solution with value $z$, we define the following metrics:
\begin{itemize}[nolistsep]
    \item \emph{Gap to best known}:
    relative cost difference from the best-known solution, computed as $(\objvalue - \bestobjvalue)/{\bestobjvalue}$;
    \item \emph{Improvement over repaired} ($\impoverrep$):
    relative cost reduction compared to the repaired solution, computed as 
    $\impoverrep(\objvalue) = (\repairedobjvalue - \objvalue)/{\repairedobjvalue}$.
\end{itemize}
In the following, we will denote by $\BLobjvalue$ and $\GNNobjvalue$ the values of the solutions produced by the baseline and the GNN-aided strategies, respectively.

\myskip

\cref{tab:reoptimization_comparison} compares the performance of these approaches on the instances of subset B, broken down by instance set (Set 1 and Set 2) and disruption type (MB for Machine Breakdown and PS for Plant Shutdown). 
Besides gap to best known and improvement over repaired, the table reports the percentage of instances for which the GNN-aided approach wins over the baseline, where column “Total” is further partitioned into cases where $\impoverrepGNN$ exceeds $\impoverrepBL$ by a margin smaller than $5\%$ ($\Delta \impoverrep < 5\%$) and cases where the improvement exceeds $5\%$ ($\Delta \impoverrep \geq 5\%$). The last three columns are defined analogously, but referring to the cases where baseline has the best performance.

\begin{table}[htbp]
    \centering
    \renewcommand{\arraystretch}{1.2}
    \scriptsize
    \resizebox{\columnwidth}{!}{
    \newcommand{\MachineBreakdown}{MB}
    \newcommand{\PlantShutdown}{PS}
    \newcommand{\pT}{}
    \newcommand{\pC}{\%}
    \begin{tabular}{ll rr rr rrr rrr}
        \toprule
        \textbf{Inst.} & 
        \textbf{Disr.} & 
        \multicolumn{2}{c}{\textbf{Gap to best \pT}} & 
        \multicolumn{2}{c}{\textbf{Imp. over rep. \pT}} & 
        \multicolumn{3}{c}{\textbf{(G) win over (B) \pT}} & 
        \multicolumn{3}{c}{\textbf{(G) loss vs (B) \pT}} \\
        \cmidrule(lr){3-4} \cmidrule(lr){5-6} \cmidrule(lr){7-9} \cmidrule(lr){10-12}
        \textbf{Set} &
        \textbf{Type} & 
        \multicolumn{1}{c}{\textbf{(B)}} & \multicolumn{1}{c}{\textbf{(G)}} & 
        $\boldsymbol{\impoverrepBL}$ & 
        $\boldsymbol{\impoverrepGNN}$ & 
        \multicolumn{1}{c}{\textbf{Total}} & 
        $\boldsymbol{\Delta \impoverrep \! < \! 5\%}$ & 
        $\boldsymbol{\Delta \impoverrep \! \geq \! 5\%}$ & 
        \multicolumn{1}{c}{\textbf{Total}} & 
        $\boldsymbol{\Delta \impoverrep \! < \! 5\%}$ & 
        $\boldsymbol{\Delta \impoverrep \! \geq \! 5\%}$ \\
        \midrule
        \multirow{3}{*}{Set 1} 
        & \MachineBreakdown & 20.93\pC & 4.56\pC & 13.87\pC & 25.29\pC & 94.17\pC & 19.83\pC & 74.33\pC &  5.83\pC &  5.00\pC & 0.83\pC \\
        &    \PlantShutdown &  8.30\pC & 2.93\pC & 10.66\pC & 15.07\pC & 82.00\pC & 36.00\pC & 46.00\pC & 18.00\pC & 16.50\pC & 1.50\pC \\
        &               All & 17.77\pC & 4.15\pC & 13.07\pC & 22.74\pC & 91.13\pC & 23.88\pC & 67.25\pC &  8.88\pC &  7.88\pC & 1.00\pC \\
        \midrule
        \multirow{3}{*}{Set 2} 
        & \MachineBreakdown & 18.25\pC & 9.61\pC & 11.99\pC & 17.84\pC & 70.46\pC & 25.44\pC & 45.02\pC & 29.54\pC & 20.64\pC & 8.90\pC \\
        &    \PlantShutdown &  6.37\pC & 3.92\pC & 10.24\pC & 12.18\pC & 63.30\pC & 41.49\pC & 21.81\pC & 36.70\pC & 30.85\pC & 5.85\pC \\
        &               All & 15.27\pC & 8.18\pC & 11.55\pC & 16.42\pC & 68.67\pC & 29.47\pC & 39.20\pC & 31.33\pC & 23.20\pC & 8.13\pC \\
        \midrule
        \multirow{3}{*}{All} 
        & \MachineBreakdown & 19.63\pC & 7.00\pC & 12.96\pC & 21.69\pC & 82.70\pC & 22.55\pC & 60.15\pC & 17.30\pC & 12.56\pC & 4.73\pC \\
        &    \PlantShutdown &  7.36\pC & 3.41\pC & 10.46\pC & 13.67\pC & 72.94\pC & 38.66\pC & 34.28\pC & 27.06\pC & 23.45\pC & 3.61\pC \\
        &               All & 16.56\pC & 6.10\pC & 12.33\pC & 19.68\pC & 80.26\pC & 26.58\pC & 53.68\pC & 19.74\pC & 15.29\pC & 4.45\pC \\
        \bottomrule
    \end{tabular}
    }
    \caption{Results of the comparison between GNN-aided (G) and baseline (B) strategies.}
    \label{tab:reoptimization_comparison}
\end{table}

The results show that the baseline approach produces solutions with a large gap to the best-known solution; the average value of this gap equals 16.56\%. 
In contrast, the GNN-aided approach substantially reduces the gap to the best-known solution: it never exceeds 9.61\% on a single group of homogeneous instances, and has an average value equal to 6.10\%.
Overall, the GNN-aided strategy outperforms the baseline for most of the instances, achieving a global average cost improvement of 19.68\% compared to 12.33\% for the baseline, within the same time limit.
Specifically, the GNN strategy wins in 80.26\% of cases. 
Notably, a significant portion of these wins (53.68\% out of 80.26\%) are substantial ($\Delta \impoverrep \geq 5\%$), while the vast majority of losses (15.29\% out of 19.74\%) are marginal ($\Delta \impoverrep < 5\%$).

Breaking down the results by instance set and disruption type reveals two key trends. First, the GNN-aided approach yields larger improvements on the uniform instances of Set 1 (on average, $\impoverrepGNN = 22.74\%$) compared to the heterogeneous instances of Set 2 (on average, $\impoverrepGNN = 16.42\%$).
This is likely because the structural homogeneity of instances in Set 1, all having 3 machines and 30 periods, facilitates the GNN in predicting good neighborhoods. Indeed, the number of instances of this type in the training Set A is proportionally larger. However, this also confirms the generalization quality of our GNN, a remarkable feature for the approach.
Second, regarding disruption types, the GNN-aided approach is more effective for machine breakdowns (on average, $\impoverrepGNN = 21.69\%$) than for plant shutdown (on average, $\impoverrepGNN = 13.67\%$).
When a machine breakdown occurs, the disruption affects a single machine, and recovery involves reallocating critical productions from the broken machine to alternative ones.
These adjustments involve a limited number of setup variables and are easily mastered by the GNN’s predictive capabilities.
In contrast, when a plant shutdown occurs, the disruption is global as productions across all the machines are canceled for a few periods. 
In this case, recovery requires an overall restructuring of the production, involving simultaneous adjustments across machines and propagating over subsequent time frames.
The improvement of the best-known solution over the repaired one is more consistent for MB instances rather than PS instances. For the latter set, although there is smaller room for improvement, the GNN-aided approach consistently outperforms the baseline, halving the gap to the best-known solution.

\subsection{Evaluating the integration of imperfect GNN predictions}

In this section, we computationally evaluate the reliability of the GNN prediction and its integration within the reoptimization framework.
Specifically, we investigate the impact of the setup selection budget $\nbselectedsetups$.
We compare the proposed configuration ($\nbselectedsetups = 30$) with two alternatives to define the subset $\setselectedsetups$ of free setup variables:
\begin{itemize}[nolistsep]
    \item \emph{Tight GNN-aided approach (T).}
    Assuming the GNN is highly accurate, the tight approach optimistically restricts the solver to the GNN's top predictions by setting the budget $\nbselectedsetups$ to exactly the maximum allowed changes $\nbsetupchanges$, i.e. $\nbselectedsetups = \nbsetupchanges = 10$.
    \item \emph{Perfect-predictor approach (P).}
    This approach is obtained by substituting the GNN predictions with the exact ground-truth labels known from the long-term solution. 
    Thus, it simulates an ideal perfect predictor capable of exactly forecasting the $\nbsetupchanges$ setup variables to change, establishing a theoretical ceiling for the tight GNN-aided fix-and-optimize strategy.
\end{itemize}
In the following, we will denote by $\Tobjvalue$ and $\Pobjvalue$ the values of the solutions produced by the tight and the perfect-predictor approaches, respectively.

\medskip

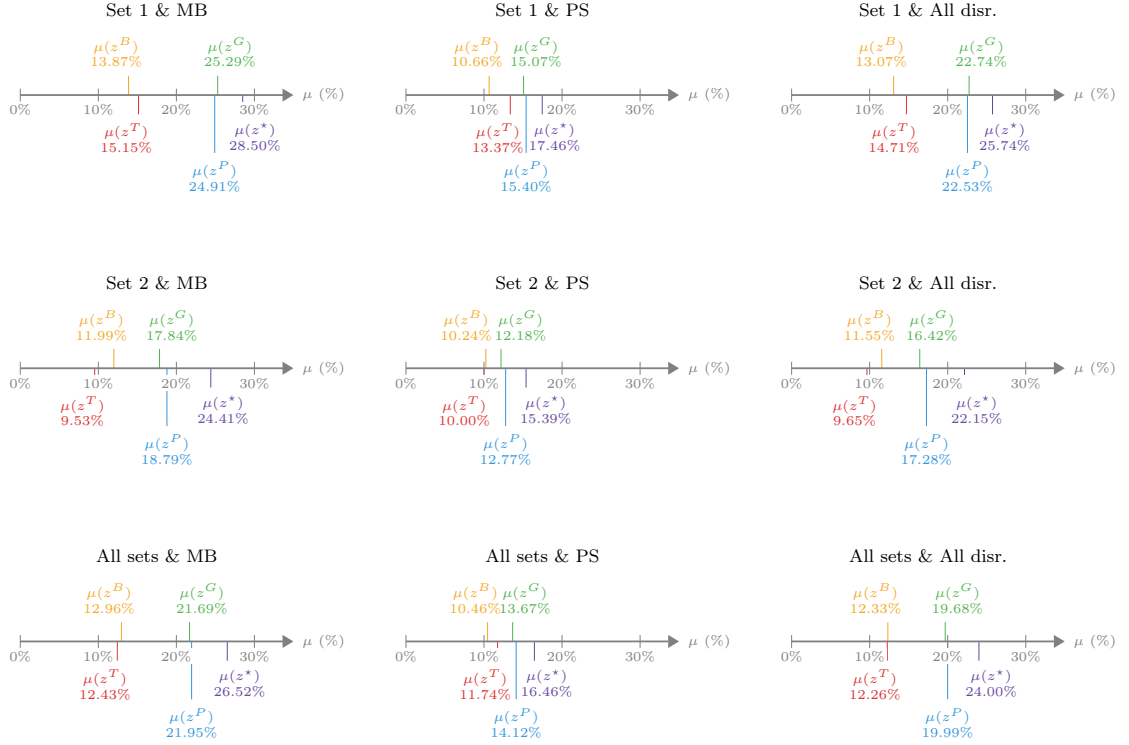
\begin{figure}[h]
    \centering
    \scalebox{.85}{
    \definecolor{colorBaseline}{RGB}{245, 166, 35}  
    \definecolor{colorGNN}{RGB}{80, 180, 80}
    \definecolor{colorLongTerm}{RGB}{100, 70, 160}
    \definecolor{colorTight}{RGB}{220, 50, 50}
    \definecolor{colorPP}{RGB}{50, 150, 220}
    
    \newcommand{\CaseTitle}[2]{{#1} \& {#2}}
    \newcommand{\MB}{MB}
    \newcommand{\PS}{PS}
    
    \begin{tikzpicture}
    
        \def\axiswidth{4.25} 
        \def\maxval{35}
        
        \newcommand{\ImprovementAxis}[9][.5em]{
            \begin{scope}[shift={(#2,#3)}]
                \node[anchor=south, font=\scriptsize] at (\axiswidth/2, 1.1) {#4};
                
                \draw[thick, gray, -Triangle] (0,0) -- (\axiswidth, 0)
                    node[right, font=\tiny, text=gray] {$\impoverrep$ (\%)};
                    
                \foreach \x in {0, 10, 20, 30} {
                    \draw[gray] ({\x/\maxval*\axiswidth}, 0.1) -- ({\x/\maxval*\axiswidth}, -0.1) 
                        node[anchor=north, font=\tiny, text=gray, fill=white, inner sep=1pt] {\x\%};
                }
                
                \begin{scope}[on background layer]
                
                    \draw[colorBaseline] ({#5/\maxval*\axiswidth}, 0) -- ({#5/\maxval*\axiswidth}, 0.3)
                        node[anchor=south, font=\tiny, align=center, xshift=-#1] {$\impoverrepBL$ \\ #5\%};
                    \draw[colorGNN] ({#7/\maxval*\axiswidth}, 0) -- ({#7/\maxval*\axiswidth}, 0.3)
                        node[anchor=south, font=\tiny, align=center, xshift=#1] {$\impoverrepGNN$ \\ #7\%};
    
                    \draw[colorPP] ({#8/\maxval*\axiswidth}, 0) -- ({#8/\maxval*\axiswidth}, -0.9)
                        node[anchor=north, font=\tiny, align=center] {$\impoverrepP$ \\ #8\%};
                    \draw[colorTight] ({#6/\maxval*\axiswidth}, 0) -- ({#6/\maxval*\axiswidth}, -0.3)
                        node[anchor=north, font=\tiny, align=center, xshift=-#1] {$\impoverrepT$ \\ #6\%};
                    \draw[colorLongTerm] ({#9/\maxval*\axiswidth}, 0) -- ({#9/\maxval*\axiswidth}, -0.3)
                        node[anchor=north, font=\tiny, align=center, xshift=#1] {$\impoverrepLG$  \\ #9\%};
                        
                \end{scope}
            \end{scope}
        }
    
        \ImprovementAxis{0}{8.5}{\CaseTitle{Set 1}{\MB}}{13.87}{15.15}{25.29}{24.91}{28.50}
        \ImprovementAxis{6}{8.5}{\CaseTitle{Set 1}{\PS}}{10.66}{13.37}{15.07}{15.40}{17.46}
        \ImprovementAxis{12}{8.5}{\CaseTitle{Set 1}{All disr.}}{13.07}{14.71}{22.74}{22.53}{25.74}
    
        \ImprovementAxis{0}{4.25}{\CaseTitle{Set 2}{\MB}}{11.99}{9.53}{17.84}{18.79}{24.41}
        \ImprovementAxis[.8em]{6}{4.25}{\CaseTitle{Set 2}{\PS}}{10.24}{10.00}{12.18}{12.77}{15.39}
        \ImprovementAxis{12}{4.25}{\CaseTitle{Set 2}{All disr.}}{11.55}{9.65}{16.42}{17.28}{22.15}
    
        \ImprovementAxis{0}{0}{\CaseTitle{All sets}{\MB}}{12.96}{12.43}{21.69}{21.95}{26.52}
        \ImprovementAxis{6}{0}{\CaseTitle{All sets}{\PS}}{10.46}{11.74}{13.67}{14.12}{16.46}
        \ImprovementAxis{12}{0}{\CaseTitle{All sets}{All disr.}}{12.33}{12.26}{19.68}{19.99}{24.00}
    
    \end{tikzpicture}
    }
    \caption{Visual comparison of average improvement over repaired $\impoverrep$ for the GNN-aided alternatives, baseline, and best-known solutions across different instance sets and disruption types.}
    \label{fig:improvements_over_repaired}
\end{figure}

\cref{fig:improvements_over_repaired} visually summarizes the average improvement $\impoverrep$ over the repaired solution  achieved by these different configurations.
For each combination of instance set and disruption type, it compares the results obtained with the proposed GNN-aided approach ($\nbselectedsetups = 30$), the tight approach, and the perfect predictor approach, alongside the baseline and long-term ones.

The relative positions of the metrics in \cref{fig:improvements_over_repaired} reveal interesting insights about the integration of the GNN predictions in a MILP reoptimization process under strict time limits. 

First, the improvement obtained with the tight approach is unsatisfactory and frequently worse than the baseline. 
Because the GNN is limited to identifying exactly $\nbsetupchanges = 10$ variables, any misclassification may then trap the solver in a restricted sub-optimal space. 
This shows that rigidly integrating imperfect GNN predictions into an optimization scheme can be counterproductive.

Conversely, our proposed approach achieves performance comparable to -- and occasionally exceeding -- what could be obtained by an ideal perfect predictor. 
By setting the budget to $\nbselectedsetups = 30$, we provide the solver with a “buffer zone”. 
Even if the GNN's predictions contain noise, the critical $\nbsetupchanges$ setup variables are very likely included within the top 30. 
Then, the MILP solver can efficiently explore the buffered $\setselectedsetups$ to identify the best subset of setup variables to change

Interestingly, in certain cases (such as Set 1 \& Machine Breakdown), the buffered GNN approach actually outperforms the perfect predictor one. 
While the latter points directly to the global long-term optimum, finding that specific solution within a tight 10wu time limit may be computationally prohibitive.
The former, by offering a slightly larger solution space, allows the solver to find high-quality alternative solutions that are much easier to reach within the short time limit.

\myskip

Ultimately, these results demonstrate that our GNN-aided approach is most effective when the GNN is used to define a limited yet flexible subset $\setselectedsetups$ rather than an optimistically tight one.
By providing a buffered neighborhood constraint, the approach successfully leverages the prediction capabilities of the GNN, while retaining enough combinatorial flexibility for the MILP solver to compensate its inaccuracies.



\section{Conclusion}
\label{sec:conclusion}

In this paper, we considered the problem of reoptimizing a given MILP after an unpredictable event made a nominal solution infeasible. 
In this setting, the requirement to compute a new solution within a very limited time budget makes solving the MILP on the disrupted instance impractical.
At the same time, applying a simple repairing heuristic often fails to provide a high-quality alternative.
We therefore introduced a novel reoptimization strategy aided by Machine Learning (ML).
More specifically, our fix-and-optimize strategy is based on very large neighborhood search scheme and incorporates a Graph Neural Network (GNN) for capturing the intricate relationships among the instance, the nominal solution, and the disruption.
We evaluated the proposed approach on a variant of the Lot Sizing Problem, under disruption scenarios in which some machines become unavailable for a certain period of time.
Computational results show that the proposed approach outperforms a baseline approach that simply solves the MILP formulation for the disrupted instance, consistently producing better solutions within the same limited time budget.
Moreover, the experiments demonstrate the ability of our approach to handle instances of different sizes, validating the choice of employing GNN within ML-aided reoptimization approaches.

\myskip

As regards future work, various directions could be explored.
As to the specific application considered in this paper, while we 
focused on machine breakdowns and plant shutdowns, the graph-based nature of our approach offers flexibility to integrate other types of disruptions (e.g., last-minute sales order cancellations, sudden demand peaks or preventative maintenance).
More in general, the generic framework of our reoptimization approach, encoding the (instance, nominal solution, disruption) triplet into a feature graph to guide a fix-and-optimize heuristic, naturally adapts to any combinatorial optimization problem.
Finally, our approach could evolve from a supervised-learning single prediction to an iterative scheme, by employing Deep Reinforcement Learning to train an agent capable of performing a number of sequential iterations. Such an agent would dynamically learn a policy to select, at each iteration, a promising neighborhood by using the improvement in the objective function as a reward signal.
The resulting approach would eliminate the reliance on computationally expensive ground-truth labels and mitigate the risks associated with single-shot prediction errors.

\section*{Acknowledgements}\label{sec:Acknowledgements}
\noindent
Enrico Malaguti and Michele Monaci were funded by the Air Force Oﬃce of Scientific Research under award number FA8655-25-1-7013.



\bibliography{references}
\pagebreak


\appendix



\section{Repairing heuristic algorithm}
\label{app:repairing_heuristic}

\vspace*{-1em}

\begin{algorithm}
    \setstretch{1}
    \caption{\rule{0pt}{4mm}Repairing heuristic \vspace{1mm}}
    \label{alg:repairing_heuristic}
    \vspace{3mm}
    \footnotesize
    \newcommand{\setperturbedmachines}{\mathcal{M}^{p}}
    \newcommand{\usedcapacity}{\capacity^{used}}
    \newcommand{\neededcapacity}{\capacity^{need}}
    \newcommand{\availablecapacity}{\capacity^{avail}}
    \newcommand{\flowbalance}{\beta}
    
    \KwIn{Original solution $\refsolution$, set of disrupted machines $\setperturbedmachines$, disruption duration $\maintenanceduration$}
    \KwOut{Repaired solution $\repairedsolution$ feasible w.r.t. disrupted instance}
    
    \BlankLine
    
    $\repairedsolution \leftarrow \refsolution$ 
    \tcp{Initialize repaired solution as a copy of original solution}
    
    \BlankLine
    
    \SetKwProg{Proc}{Procedure}{}{}\SetKwFunction{procCall}{Call}
    \Proc{CancelDisruptedProductions}{}{
        \vspace{-1mm}
        \For{$\machineindex \in \setperturbedmachines$}{
            \For{$\periodindex \in \{1, \ \dots, \ \maintenanceduration \}$}{
                \For{$\itemindex \in \{1, \ \dots, \ \nbitems\}$}{
                    $\repairedvarsetup_{\threegenericindices} \leftarrow 0$, \,
                    $\repairedvarcarryover_{\threegenericindices} \leftarrow 0$, \,
                    $\repairedvarquantity_{\threegenericindices} \leftarrow 0$ \\
                }
            }
        }
    }
    
    \BlankLine
    
    \Proc{AssessBrokenCarryOvers}{}{
        \vspace{-1mm}
        $\periodindex \leftarrow \maintenanceduration$, \,
        $\periodindex' \leftarrow \maintenanceduration + 1$ \\
        \For{$\machineindex \in \setperturbedmachines$}{
            $\usedcapacity_{\machineindex \periodindex'} \leftarrow \sum_{\itemindex=1}^{\nbitems} (\setuptime_{\itemindex} \refvarsetup_{\threegenericindices'} + \productiontime_{\itemindex} \refvarquantity_{\threegenericindices'})$ 
            \tcp{Capacity of machine $\machineindex$ used}
            \For{$\itemindex \in \{1, \ \dots, \ \nbitems\}$}{
                \tcp{If production of item $\itemindex$ at period $\periodindex'$ in $\refsolution$ depended on a now broken carry-over}
                \If{$\refvarcarryover_{\threegenericindices} = 1$ \textbf{and} $\refvarquantity_{\threegenericindices'} > 0$}{
                    $\neededcapacity \leftarrow \setuptime_{\itemindex} + \productiontime_{\itemindex} \minproduction_{\itemindex}$ 
                    \tcp{Capacity needed to setup and produce min quantity of $\itemindex$ at $\periodindex'$}
                    $\availablecapacity \leftarrow \capacity_{\machineindex \periodindex'} - (\usedcapacity_{\machineindex \periodindex'} - (\setuptime_{\itemindex} \refvarsetup_{\threegenericindices'} + \productiontime_{\itemindex} \refvarquantity_{\threegenericindices'}))$ 
                    \tcp{Capacity available on $\machineindex$ at $\periodindex'$}
                    \tcp{If needed capacity larger than available one, then cancel prod. of $\itemindex$ at $\periodindex'$}
                    \If{$\neededcapacity > \availablecapacity$}{
                        $\repairedvarsetup_{\threegenericindices'} \leftarrow 0$, \,
                        $\repairedvarquantity_{\threegenericindices'} \leftarrow 0$ \\
                    }
                    \tcp{Else, produce min quantity of $\itemindex$ at $\periodindex'$}
                    \Else{
                        $\repairedvarsetup_{\threegenericindices'} \leftarrow 1$, \,
                        $\repairedvarquantity_{\threegenericindices'} \leftarrow \minproduction_{\itemindex}$ \\
                    }
                }
            }
        }
    }
    
    \BlankLine
    
    \Proc{RecomputeInventoryFlow}{}{
        \vspace{-1mm}
        \For{$\itemindex \in \{1, \ \dots, \ \nbitems\}$}{
            $\repairedvarinventory_{\itemindex 0} \leftarrow \refvarinventory_{\itemindex 0}$ \\
            \For{$t \in \{1, \ \dots, \ \nbperiods\}$}{
                $\flowbalance \leftarrow \repairedvarinventory_{\itemindex (\periodindex - 1)} + \big(\sum_{\machineindex = 1}^{\nbmachines} \repairedvarquantity_{\threegenericindices} \big) - \demand_{\itemindex \periodindex}$ 
                \tcp{Flow balance of item $\itemindex$ at period $\periodindex$}
                \If{$\flowbalance \ge 0$}{
                    $\repairedvarinventory_{\itemindex \periodindex} \leftarrow \flowbalance$, \,
                    $\repairedvarlostsales_{\itemindex \periodindex} \leftarrow 0$ \\
                }
                \Else{
                    $\repairedvarinventory_{\itemindex \periodindex} \leftarrow 0$, \,
                    $\repairedvarlostsales_{\itemindex \periodindex} \leftarrow \min(-\flowbalance, \ \demand_{\itemindex \periodindex})$ \\
                }
            }
        }
    }
    
    \BlankLine
    
    \procCall{CancelDisruptedProductions} \\
    \procCall{AssessBrokenCarryOvers} \\
    \procCall{RecomputeInventoryFlow} \\
    
    \BlankLine
    
    \KwRet{$\repairedsolution$}
    \vspace{2mm}
\end{algorithm}


\pagebreak



\section{Notations}

\begin{table}[ht!]
    \centering
    \renewcommand{\arraystretch}{1.2}
    \footnotesize
    \begin{tabular}{c c l}
        \toprule
        \textbf{Category} & \textbf{Symbol} & \textbf{Description} \\ 
        \midrule
        \multirow{3}{*}{Indices} 
        & $\itemindex$ & Item index \\
        & $\machineindex$ & Machine index \\
        & $\periodindex$ & Period index \\ 
        \midrule
        \multirow{12}{*}{Instance} 
        & $\nbmachines$ & Number of machines \\
        & $\nbitems$ & Number of items \\
        & $\nbperiods$ & Number of periods \\ 
        \cmidrule(l){2-3}
        & $\setupcost_{\itemindex}$ & Setup cost of item $\itemindex$ \\
        & $\productioncost_{\itemindex}$ & Unit production cost of item $\itemindex$ \\
        & $\inventorycost_{\itemindex}$ & Unit inventory cost of item $\itemindex$ \\
        & $\lostsalescost_{\itemindex}$ & Unit lost sales cost of item $\itemindex$ \\
        & $\demand_{\itemindex \periodindex}$ 
        & Demand of item $\itemindex$ at period $\periodindex$ \\
        & $\setuptime_{\itemindex}$ & Setup time of item $\itemindex$ \\
        & $\productiontime_{\itemindex}$ & Unit production time of item $\itemindex$ \\
        & $\capacity_{\machineindex \periodindex}$ 
        & Capacity of machine $\machineindex$ at period $\periodindex$ \\
        & $\minproduction_{\itemindex}$ & Minimum production of item $\itemindex$ \\ 
        \midrule
        \multirow{1}{*}{Disruption}
        & $\maintenanceduration$ & Disruption duration \\
        \midrule
        \multirow{2}{*}{Neighborhood}
        & $\shorttermhorizon$ & Number of periods in the short-term horizon \\
        & $\nbsetupchanges$ & Bound on the number of setup changes \\ 
        \midrule
        \multirow{9}{*}{Solution}
        & $\varsetup_{\threegenericindices}$ 
        & Setup decision variable \\
        & $\varcarryover_{\threegenericindices}$ 
        & Carry-over decision variable \\
        & $\varquantity_{\threegenericindices}$ 
        & Production quantity decision variable \\
        & $\varinventory_{\itemindex \periodindex}$ 
        & Inventory decision variable \\
        & $\varlostsales_{\itemindex \periodindex}$ 
        & Lost-sales decision variable \\ 
        \cmidrule(l){2-3}
        & $\refsolution$ 
        & Original solution (feasible w.r.t. original optimization problem) \\
        & $\repairedsolution$ 
        & Repaired solution (feasible w.r.t. reoptimization problem) \\ 
        & $\newsolution$ 
        & New solution (feasible w.r.t. reoptimization problem) \\ 
        & $\optimalsolution$
        & Optimal (or near-optimal) new solution (feasible w.r.t. reoptimization problem) \\
        \midrule
        \multirow{10}{*}{Prediction}
        & $\embeddingdimension$ 
        & GNN embedding dimension \\
        & $\nbconvolutionblocks$ 
        & GNN number of convolution blocks \\ 
        \cmidrule(l){2-3}
        & $\groundtruthlabel_{\machineitemperiodindices}$
        & Ground-truth label associated to setup variable $\varsetup_{\threegenericindices}$ \\
        & $\score_{\threegenericindices}$ 
        & Change-likelihood score associated to setup variable $\varsetup_{\threegenericindices}$ \\
        & $\nbselectedsetups$ 
        & Setup selection budget \\
        & $\setselectedsetups$ 
        & Subset of selected free setup decision variables in the short-term horizon \\ 
        \cmidrule(l){2-3}
        & $\alpha$ 
        & Focal loss class-balancing parameter \\ 
        & $\gamma$ 
        & Focal loss focusing parameter \\ 
        \cmidrule(l){2-3}
        & $\learningrate$ 
        & Adam optimizer initial learning rate  \\
        \bottomrule
    \end{tabular}
    \caption{Table of notations}
    \label{tab:notations}
\end{table}

\end{document}